\documentclass[aos,citesort,dvips]{arximspdf}

\doi{10.1214/10-AOS809}
\volume{38}
\issue{5}
\pubyear{2010}
\firstpage{2884}
\lastpage{2915}

\makeatletter

\newcommand{\XW}{\mathit{XW}}

\newtheorem{theorem}{Theorem}
\newtheorem{corollary}{Corollary}
\newtheorem{lemma}{Lemma}
\newtheorem{proposition}{Proposition}
\newproclaim{remark}{Remark}

\makeatother

\begin{document}
\begin{frontmatter}

\title{Bootstrap consistency for general semiparametric $M$-estimation}
\runtitle{Semiparametric bootstrap consistency}

\begin{aug}
\author[A]{\fnms{Guang} \snm{Cheng}\corref{}\thanksref{t1}\ead[label=e1]{chengg@purdue.edu}} and
\author[B]{\fnms{Jianhua Z.} \snm{Huang}\thanksref{t2}\ead[label=e2]{jianhua@stat.tamu.edu}}
\runauthor{G. Cheng and J. Z. Huang}
\affiliation{Purdue University and Texas A\&M University}
\address[A]{Department of Statistics\\
Purdue University\\
West Lafayette, Indiana 47907-2066\\
USA\\
\printead{e1}}
\address[B]{Department of Statistics\\
Texas A\&M University\\
College Station, Texas 77843-3143\\
USA\\
\printead{e2}}
\end{aug}

\thankstext{t1}{Supported by NSF Grant DMS-09-06497.}
\thankstext{t2}{Supported in part by NSF Grants DMS-06-06580, DMS-09-07170,
NCI Grant CA57030 and Award Number KUS-CI-016-04, made by King Abdullah
University of Science and Technology (KAUST).}

% HISTORY:
\received{\smonth{10} \syear{2009}}

% ABSTRACT
%
\begin{abstract}
Consider $M$-estimation in a semiparametric model that is\break characterized
by a Euclidean parameter of interest and an infinite-dimensional
nuisance parameter. As a general purpose approach to statistical
inferences, the bootstrap has found wide applications in semiparametric
$M$-estimation and, because of its simplicity, provides an attractive
alternative to the inference approach based on the asymptotic
distribution theory. The purpose of this paper is to provide
theoretical justifications for the use of bootstrap as a semiparametric
inferential tool. We show that, under general conditions, the bootstrap
is asymptotically consistent in estimating the distribution of the
$M$-estimate of Euclidean parameter; that is, the bootstrap distribution
asymptotically imitates the distribution of the $M$-estimate. We also
show that the bootstrap confidence set has the asymptotically correct
coverage probability. These general conclusions hold, in particular,
when the nuisance parameter is not estimable at root-$n$ rate, and apply
to a broad class of bootstrap methods with exchangeable bootstrap
weights. This paper provides a first general theoretical study of the
bootstrap in semiparametric models.
\end{abstract}

% KEYWORDS
%

\end{frontmatter}

%s1 ###
\section{Introduction}
Due to its flexibility, semiparametric modeling has provided a
powerful statistical modeling framework for complex data, and
proven to be useful in a variety of contexts, see
\cite{bmm08,h99,zl07,c08a,zy08}. Semiparametric models are indexed
by a Euclidean parameter of interest
$\theta\in\Theta\subset\mathbb{R}^d$ and an infinite-dimensional
nuisance function $\eta$ belonging to a Banach space $\mathcal{H}$
with a norm~\mbox{$\|\cdot\|$}. $M$-estimation, including the maximum
likelihood estimation as a special case, refers to a general
method of estimation, where the estimates are obtained by
optimizing some objective functions \cite{dhp06,wz07,mk05a}.
The asymptotic theories and inference procedures for
semiparametric maximum likelihood estimation, or more generally
$M$-estimation, have been extensively studied in
\cite{bkrw98,k08a,mv00,mk05a,lkf05,dkl05}.

It is well known that the asymptotic inferences of semiparametric
models often face practical challenges. In particular, the
confidence set construction and the asymptotic variance estimation
of the estimator for the Euclidean parameter both involve
estimating and inverting a hard-to-estimate infinite-dimensional
operator. The difficulty in dealing with such an
infinite-dimensional operator motivated the development of the
profile sampler \cite{lkf05,ck08a,ck08b}, where the inference of
the Euclidean parameter is based on sampling from the posterior of
the profile likelihood \cite{lkf05}. However, because of the way
it is designed, the profile sampler method has the typical caveats
of the Bayesian methods. First, one needs to specify a prior
distribution. Second, since the Markov chain Monte Carlo (McMC) is
used for sampling from the posterior distribution, there are a
number of controversial issues in generating the stationary Markov
chain. For example, it is considerably difficult to determine the
burn-in period and stopping time of the chain \cite{gcsr03}. In
particular, it may take a long time to run the Markov chain in
order to give accurate inferences for $\theta$ when $\eta$ is
estimable at a slow convergence rate \cite{ck08a,ck08b}.
Moreover, when the sample size is small, the profile likelihood
may become nonsmooth or may not approximate well the desired
parabolic form, violating the main theoretical basis of the
profile sampler.

On the other hand, as a general data-resampling based statistical
inference tool, the bootstrap method does not have the drawbacks
of the profile sampler; see \cite{cp09,mk05a,s06,hhms04,klf04,yjs08}
for its application in semiparametric models. In fact, the
bootstrap method has several methodological advantages over the
profile sampler: it is straightforward to implement; there is no
need to specify a prior distribution and to check Markov chain
convergence. In addition, the finite sample performance of the
bootstrap can be controlled by choosing from a rich pool of
resampling techniques; see Section 3 of \cite{pw93}. Moreover,
unlike the profile sampler which focuses on $\theta$, one can make
bootstrap inferences for both $\theta$ and $\eta$.

Unfortunately, a systematic theoretical study on the bootstrap
inference in semiparametric models is almost nonexistent,
especially when the nuisance function parameter $\eta$ is not
$\sqrt{n}$ estimable, despite the rich literature on the bootstrap
theory for parametric models \cite{bf81,s81,h92,mn92}. The
current literature only considered the bootstrap consistency for
the joint estimator of $(\theta,\eta)$ in some special case of
semiparametric models where $\eta$ is $\sqrt{n}$-estimable, that is,
\cite{klf04}. In a recent monograph, Kosorok pointed out that
``convergence rate and asymptotic normality results are quite
difficult to establish for the nonparametric bootstrap (based on
multinomial weights), especially for models with parameters not
estimable at the $\sqrt{n}$ rate'' \cite{k08a}. In fact, the lack of
theoretical justifications of the bootstrap in the semiparametric
context is one of the main motivations for developing the profile
sampler. The purpose of this paper is to develop a general theory
on bootstrap consistency in semiparametric models, for a broad
class of bootstrap methods including Efron's (nonparametric)
bootstrap as a special case. We focus on the inference of the
Euclidean parameter and leave study of the bootstrap inference of
the nuisance parameter for future research, although we give some
useful convergence rate results (see Section \ref{secrate}).

Our main results are summarized as follows. The semiparametric $M$-estimator
$(\widehat{\theta},\widehat{\eta})$ and the bootstrap $M$-estimator
$(\widehat{\theta}^{\ast},\widehat{\eta}^{\ast})$ are obtained by
optimizing the objective function $m(\theta,\eta)$ based on the
i.i.d. observations $(X_1,\ldots,X_n)$ and the bootstrap sample
$(X_1^{\ast},\ldots,X_n^{\ast})$, respectively:
%
%e2 ###
%e1 ###
%
\begin{eqnarray}
\label{mestimate}
(\widehat{\theta},\widehat{\eta})
&=&\mathop{\arg\sup}_{\theta\in\Theta,\eta\in\mathcal{H}}\sum_{i=1}^{n}
m(\theta,\eta)(X_{i}),\\
\label{best}
(\widehat{\theta}^{\ast},\widehat{\eta}^{\ast})
&=&\mathop{\arg\sup}_{\theta\in\Theta,\eta\in\mathcal{H}}\sum_{i=1}^{n}
m(\theta,\eta)(X_{i}^{\ast}),
\end{eqnarray}
where $(X_1^{\ast},\ldots,X_n^{\ast})$ are independent draws with
replacement from the original sample. Note that we can express
%
%e3 ###
%
\begin{equation}\label{inview}
(\widehat{\theta}^{\ast},\widehat{\eta}{}^{\ast})
=\mathop{\arg\sup}_{\theta\in\Theta,\eta\in\mathcal{H}}\sum_{i=1}^{n}W_{ni}
m(\theta,\eta)(X_{i}),
\end{equation}
and the bootstrap weights $(W_{n1},\ldots,W_{nn})\sim
\operatorname{Multinomial}(n,(n^{-1},\ldots,n^{-1}))$. In this paper, we
consider the more general \textit{exchangeable bootstrap weighting
scheme} that includes Efron's bootstrap and its smooth alternative
\cite{l93}, for example, \textit{Bayesian bootstrap}, as special cases. The
general resampling scheme was first proposed in \cite{r81}, and
extensively studied by \cite{bb95}, who suggested the name
``\textit{weighted bootstrap},'' and in \cite{pw93,mn92}. Note that
other variations of Efron's bootstrap are also studied in
\cite{cb05} using the term ``\textit{generalized bootstrap}.'' The
practical usefulness of the more general scheme is well-documented
in the literature. For example, in semiparametric survival models,
for example, Cox regression model, the nonparametric bootstrap often
gives many ties when it is applied to censored survival data due
to its ``discreteness'' and the general weighting scheme comes to
the rescue. As one main contribution of the paper, we show that
the nonparametric bootstrap distribution of
$\sqrt{n}(\widehat{\theta}^{\ast}-\widehat{\theta})$, conditional
on the observed data, asymptotically imitates the distribution of
$\sqrt{n}(\widehat{\theta}-\theta_0)$, where $\theta_0$ is the
true value of $\theta$. As a consequence, we also establish the
consistency of the bootstrap confidence set of $\theta$, which
means that the coverage probability converges to the nominal
level. Our results hold when the estimate of the nuisance function
has either root-$n$ or slower than root-$n$ convergence rate. This
paper can also be viewed as a nontrivial extension of \cite{cb05}
to account for the presence of an infinite-dimensional nuisance
parameter.

In a related paper, Ma and Kosorok \cite{mk05a} obtained some
theoretical results when the bootstrap weights are assumed to be
i.i.d. There is a crucial difference between their work and ours:
They treated the bootstrap estimator as the regular weighted
estimator and used the unconditional arguments rather than the
usual conditional arguments as we employ in this paper. Note that
the i.i.d. assumption rules out all interesting bootstrap schemes
considered in this paper, and their theoretical approach cannot be
extended to obtain our results. Indeed, they stated in the paper
that the independence assumption makes their proofs easier and the
relaxation to the dependent weights appears to be quite difficult.
Another related work is the piggyback bootstrap \cite{dkl05},
which is invented solely to draw inferences for the functional
parameter $\eta$ when it is $\sqrt{n}$-estimable. The piggyback
bootstrap is not the standard bootstrap and relies on a valid
random draw from the asymptotic distribution of the estimate of
$\theta$, which is hard to estimate in general. Other related work
includes interesting results on bootstrap (in)-consistency in
nonparametric estimation; see \cite{wz96,k08b,sbw09}. An $m$
out of $n$ bootstrap was developed for nonstandard $M$-estimation
with nuisance parameters in parametric models \cite{lp06}.
%For parametric models, an $m$ out of $n$ bootstrap was developed
%when the nuisance parameter is not $\sqrt{n}$-estimable

%The rest of the paper is organized as follows.
Section \ref{background} provides
the necessary background of $M$-estimation in semiparametric models.
Our main results, including the bootstrap consistency theorem, are
presented in Section \ref{bootcons}. Sections \ref{veri-con}
and \ref{secrate} discuss how to verify various technical
conditions needed for the main results. Section \ref{example}
illustrates the applications of our main results in three
examples. Section \ref{pfthm2} contains the proof of the main
results in Section \ref{bootcons}. Some useful lemmas and
additional proofs are postponed to \hyperref[app]{Appendix}.

%

%s2 ###
\section{Background}\label{background}
We first introduce a paradigm for the semiparametric $M$-estimation
\cite{mk05a,wz07}, which parallels the efficient
influence function paradigm used for the MLEs [where $m(\theta,\eta)$
is the log likelihood].
Next, we present the model assumptions needed
for the remainder of the paper, and, finally, we review some known
results on the asymptotic distribution of semiparametric
$M$-estimators, which are needed in studying the asymptotic
properties of the bootstrap.

Let
\[
m_{1}(\theta,\eta)=\frac{\partial}{\partial\theta}m(\theta,\eta)
\quad\mbox{and}\quad m_{2}(\theta,\eta)[h]=
\frac{\partial}{\partial t}m(\theta,\eta(t)) \Big|_{t=0},
\]
where $h$ is\vspace*{1pt} a ``direction'' along which $\eta(t)\in\mathcal{H}$
approaches $\eta$ as $t\rightarrow0$, running through some index
set $\mathbf{H}\subseteq L_{2}^{0}(P_{\theta,\eta})$. Similarly, we
also define
\begin{eqnarray*}
m_{11}(\theta,\eta)&=&\frac{\partial}{\partial\theta}m_{1}
(\theta,\eta) \quad\mbox{and}\quad m_{12}(\theta,\eta)[h]=
\frac{\partial}{\partial t} m_{1}(\theta,\eta(t)) \Big|_{t=0},\\
m_{21}(\theta,\eta)[h]&=&\frac{\partial}{\partial\theta
}m_{2}(\theta,
\eta)[h] \quad\mbox{and}\quad m_{22}(\theta,\eta)[h, g]=
\frac{\partial}{\partial t}m_{2}(\theta,\eta_2(t))[h] \Big|_{t=0},
\end{eqnarray*}
where $h,g\in\mathbf{H}$ and $(\partial/\partial
t)\eta_{2}(t) |_{t=0}=g$. Define
\begin{eqnarray*}
m_{2}(\theta,\eta)[H]&=& (m_{2}(\theta,\eta)
[h_{1}],\ldots,m_{2}(\theta,\eta)
[h_{d}] )',\\
m_{22}[H,h]&=&(m_{22}(\theta,\eta)[h_1,h],\ldots,m_{22}
(\theta,\eta)[h_d,h])',
\end{eqnarray*}
where $H=(h_1,\ldots,h_d)$ and $h_{j}\in\mathbf{H}$ for
$j=1,\ldots,d$. Assume there exists an
\[
H^{\dag}(\theta,\eta)=(h_1^{\dag}(\theta,\eta),\ldots,h_d^{\dag}
(\theta,\eta))',
\]
where each
$h_{j}^{\dag}(\theta,\eta)\in\mathbf{H}$, such that for any
$h\in\mathbf{H}$
%
%e4 ###
%
\begin{equation}\label{lfseq}
E_{\theta,\eta} \{m_{12}(\theta,\eta)[h]-
m_{22}(\theta,\eta)[H^{\dag},h] \}=0.
\end{equation}

Following the idea of the \textit{efficient score function}, we
define the function
\[
\widetilde{m}(\theta,\eta)=
m_{1}(\theta,\eta)-m_{2}(\theta,\eta) [H^{\dag}(\theta,\eta)].
\]
We assume that the observed data are from the probability space
$(\mathcal{X}, \mathcal{A}, P_X)$, and that
%
%e5 ###
%
\begin{equation}\label{zeroeq}
P_X\widetilde{m}(\theta_0,\eta_0)=0,
\end{equation}
where $P_Xf$ is the customary operator notation defined as $\int
f\,dP_X$. The assumption (\ref{zeroeq}) is common in semiparametric
$M$-estimation \cite{wz07,mk05a} and usually holds by the model
specifications, for example, the semiparametric regression models with
``panel count data'' \cite{wz07}. In particular, when
$m(\theta,\eta)=\log \operatorname{lik}(\theta,\eta)$, (\ref{zeroeq}) trivially
holds and $\widetilde{m}(\theta,\eta)$ becomes the well studied
\textit{efficient score function} for $\theta$ in semiparametric
models, see \cite{bkrw98}. Since
$(\widehat{\theta},\widehat{\eta})$ is assumed to be the maximizer
of $\sum_{i=1}^{n}m(\theta,\eta)(X_i)$,
$(\widehat{\theta},\widehat{\eta})$ satisfies
%
%e6 ###
%
\begin{equation}\label{exactsol0}
\mathbb{P}_{n}\widetilde{m}(\widehat{\theta},\widehat{\eta
})=0,
\end{equation}
where $\mathbb{P}_{n}f$ denotes $\sum_{i=1}^{n}f(X_i)/n$. The
theory developed in this paper is general enough to deal with the
case that $(\widehat{\theta},\widehat{\eta})$ is not the exact maximizer.
Instead of (\ref{exactsol0}), we only assume the
following ``nearly-maximizing'' condition
%
%e7 ###
%
\begin{equation} \label{exactsol}
\mathbb{P}_{n}\widetilde{m}(\widehat{\theta},\widehat{\eta})=
o_{P_X}^o(n^{-1/2}),
\end{equation}
where the superscript ``$o$'' denotes the outer probability.

Throughout the rest of the paper, we use the shortened notation
$H_0^{\dag}=H^{\dag}(\theta_0,\eta_0)$,
$\widetilde{m}_{0}=\widetilde{m}(\theta_{0},\eta_{0})$ and
$\widehat{m}=\widetilde{m}(\widehat{\theta},\widehat{\eta})$. For
a probability space $(\Omega,\mathcal{A},P)$ and a map
$T\dvtx\Omega\mapsto\bar{\mathbb{R}}$ that need not be measurable, the
notation $E^{o}T$, $O_{P}^{o}(1)$, and $o_{P}^{o}(1)$ represent
the outer expectation of $T$ w.r.t. $P$, bounded and converging to
zero in outer probability,\vspace*{1pt} respectively. More precise definitions
can be found on page 6 of \cite{vw96}. Let $V^{\otimes2}$
represent $VV'$ for any vector $V$. Define $x\vee y$ ($x\wedge y$)
to be the maximum (minimum) value of $x$ and $y$.

We now state some general conditions that will be used throughout
the whole paper. We assume that the true value $\theta_{0}$ of the
Euclidean parameter is an interior point of the compact set
$\Theta$. Define
%
%e9 ###
%e8 ###
%
\begin{eqnarray}
\label{a}
A&=&P_X\{(\partial/\partial\theta)|_{\theta=\theta_0}\widetilde{m}
(\theta,\eta_0)\}=P_X\{m_{11}(\theta_{0},\eta_{0})-
m_{21}(\theta_{0}, \eta_{0})[H_0^{\dag}]\},\\
\label{b}
B&=&\operatorname{Var}\{\widetilde{m}_0(X)\}=P_X\bigl[\{m_{1}(\theta_{0},\eta_{0})-m_{2}
(\theta_{0},\eta_{0})[H_0^{\dag}]\}^{\otimes2}\bigr].
\end{eqnarray}

\begin{enumerate}[I.]
\item[I.] \textit{Positive information condition}: the matrices $A$
and $B$ are both nonsingular.
\end{enumerate}

Condition I above is used to ensure the nonsingularity of the
asymptotic variance of $\widehat{\theta}$, which will be shown to
be $A^{-1}B (A^{-1})'$; see Proposition \ref{asythm}.

For the empirical process
$\mathbb{G}_n=\sqrt{n}(\mathbb{P}_n-P_X)$, denote its norm with
respect to a function class $\mathcal{F}_{n}$ as
$\|\mathbb{G}_{n}\|_{\mathcal{F}_{n}}={\sup_{f\in\mathcal{F}_n}}|
\mathbb{G}_n f|$. For any fixed $\delta_n>0$, define a class of
functions $\mathcal{S}_{n}$ as
%
%e10 ###
%
\begin{equation}\label{sndfn}
\mathcal{S}_{n} \equiv\mathcal{S}_n(\delta_n) =
\biggl\{\frac{\widetilde{m}(\theta_{0},\eta)-\widetilde{m}
(\theta_{0},\eta_{0})}{\|\eta- \eta_{0}\|}\dvtx
\|\eta-\eta_{0}\|\leq\delta_{n} \biggr\}
\end{equation}
and a shrinking neighborhood of $(\theta_0, \eta_0)$ as
%
%e11 ###
%
\begin{equation}\label{eq:par-neighbor}
\mathcal{C}_{n}\equiv\mathcal{C}_n(\delta_n) =
\{(\theta,\eta)\dvtx\|\theta-\theta_{0}\|\leq\delta_{n},\|\eta-\eta_{0}
\|\leq\delta_{n}\}.
\end{equation}
The next two conditions S1 and S2 imply that the empirical processes
indexed by $\widetilde{m}(\theta,\eta)$ are well behaved and
$\widetilde{m}(\theta,\eta)$ is smooth enough around
$(\theta_0,\eta_0)$.

\begin{enumerate}[S1.]
\item[S1.] \textit{Stochastic equicontinuity condition}: for any
$\delta_n\rightarrow0$,
%
%e12 ###
%
\begin{equation}\label{conmod}
\|\mathbb{G}_n\|_{\mathcal{S}_n}=O_{P_X}^{o}(1)
\end{equation}
and
%
%e13 ###
%
\begin{equation}\label{conmod2}
\mathbb{G}_{n}\bigl(\widetilde{m}(\theta,\eta)-\widetilde
{m}(\theta_{0},\eta)\bigr)
=O_{P_X}^{o}(\|\theta-\theta_{0}\|) \qquad\mbox{for }
(\theta,\eta)\in\mathcal{C}_{n}.
\end{equation}
\item[S2.] \textit{Smoothness condition}:
%
%e14 ###
%
\begin{equation}\label{smooth}
P_X\bigl(\widetilde{m}(\theta,\eta)-\widetilde{m}_0\bigr)=A
(\theta-\theta_{0})+O(\|\theta-\theta_{0}\|^{2}\vee\|\eta-\eta_{0}
\|^{2})
\end{equation}
for $(\theta,\eta)$ in some neighborhood of
$(\theta_{0},\eta_{0})$.
\end{enumerate}

For any fixed $\theta$, define
\[
\widehat{\eta}_{\theta}=\mathop{\arg\sup}_{\eta\in\mathcal{H}}\mathbb{P}_n
m(\theta,\eta).
\]
The next condition says that $\widehat{\eta}_{\theta}$ should be
close to $\eta_0$ if $\theta$ is close to $\theta_0$.

\begin{enumerate}[S3.]
\item[S3.] \textit{Convergence rate condition}: there exists a $\gamma
\in(1/4, 1/2]$ such that
%
%e15 ###
%
\begin{equation}\label{conrate}
\|\widehat{\eta}_{\widetilde{\theta}}-\eta_{0}\|=
O_{P_X}^{o}(\|\widetilde{\theta}-\theta_{0}\|\vee
n^{-\gamma})
\end{equation}
for any consistent $\widetilde{\theta}$.
\end{enumerate}

The above range requirement of $\gamma$ is always satisfied for
regular semiparametric models; see Section 3.4 of \cite{vw96}.
Verifications of conditions S1--S3 will be discussed in
Sections \ref{veri-con} and \ref{secrate}, and illustrated with
examples in Section \ref{example}.

The following proposition summarizes a known result on the
the asymptotic normality of the semiparametric $M$-estimator
$\widehat{\theta}$ \cite{k08a,mk05a,wz07}, which plays
an important role in proving bootstrap consistency in next section.
\begin{proposition}\label{asythm}
Suppose that conditions \textup{I, S1--S3} hold and that
$(\widehat{\theta},\widehat{\eta})$ satisfies (\ref{exactsol}). If
$\widehat{\theta}$ is consistent, then
%
%e16 ###
%
\begin{equation}\label{asydis-nb}
\sqrt{n}(\widehat{\theta}-\theta_{0})=-\sqrt{n}A^{-1}
\mathbb{P}_{n}\widetilde{m}_{0}+o_{P_X}^{o} (1).
\end{equation}
Consequently,
%
%e17 ###
%
\begin{equation}\label{asynor}
\sqrt{n}(\widehat{\theta}-\theta_0)\stackrel{d}{\longrightarrow}N(0,
\Sigma),
\end{equation}
where $\Sigma\equiv A^{-1}B (A^{-1})'$, $A$ and $B$ are given in
(\ref{a}) and (\ref{b}), respectively.
\end{proposition}

We assume consistency of $\widehat\theta$ in
Proposition \ref{asythm}. The consistency can usually be
guaranteed under the following ``well-separated'' condition
%
%e18 ###
%
\begin{equation}\label{sepcon}
P_Xm(\theta_{0},\eta_{0})>\sup_{(\theta,\eta)\notin G}P_X
m(\theta,\eta)
\end{equation}
for any open set $G\subset\Theta\times\mathcal{H}$ containing
$(\theta_{0},\eta_{0})$, see Theorem 5.7 in \cite{v98}. For
maximum likelihood estimation, that is, $m(\theta,\eta)=\log
\operatorname{lik}(\theta,\eta)$, it is easy to see that $A=-B$ and $\Sigma=
B^{-1}$, and thus $\Sigma^{-1}$ becomes the \textit{efficient
information matrix}.
\begin{remark}
Given any consistent estimator $\widehat{\Sigma}$ of $\Sigma$, we
have
%
%e19 ###
%
\begin{equation}\label{asynorest}
\sqrt{n}\widehat{\Sigma}^{-1/2}(\widehat{\theta}-\theta
_0)\stackrel{d}
{\longrightarrow}N(0, I)
\end{equation}
by Proposition \ref{asythm} and Slutsky's theorem. In practice,
a consistent $\widehat{\Sigma}$ can be obtained via either the
observed profile information approach \cite{mv99} or the profile
sampler approach \cite{lkf05}.
\end{remark}

%s3 ###
\section{Main results: Bootstrap consistency}\label{bootcons}
In this section, we establish the consistency of bootstrapping
$\theta$ under general conditions in the framework of
semiparametric $M$-estimation. Define
\[
\mathbb{P}_n^{\ast}f=(1/n)\sum_{i=1}^{n}W_{ni}f(X_i),
\]
where
$W_{ni}$'s are the bootstrap weights defined on the probability
space $(\mathcal{W},\Omega, P_{W})$. In view of (\ref{inview}), the
bootstrap estimator can be rewritten as
%
%e20 ###
%
\begin{equation}\label{boodef}
(\widehat{\theta}^{\ast},\widehat{\eta}^{\ast})
=\mathop{\arg\sup}_{\theta\in\Theta,\eta\in\mathcal{H}}\mathbb
{P}_n^{\ast}
m(\theta,\eta).
\end{equation}
The definition of $(\widehat{\theta}^\ast,\widehat{\eta}^\ast)$,
that is, (\ref{boodef}), implies that
%
%e21 ###
%
\begin{equation}\label{exactsol-b0}
\mathbb{P}_{n}^{\ast}\widetilde{m}(\widehat{\theta}^{\ast},
\widehat{\eta}^{\ast})=0.
\end{equation}
Similar to (\ref{exactsol}), we weaken
(\ref{exactsol-b0}) to the following ``nearly-maximizing''
condition
%
%e22 ###
%
\begin{equation}\label{exactsol-b}
\mathbb{P}_{n}^{\ast}\widetilde{m}(\widehat{\theta}^{\ast},
\widehat{\eta}^{\ast})=o_{P_{\XW}}^{o}(n^{-1/2}),
\end{equation}
where $P_{\XW}$ is a probability measure on a product space that
we will formally define later.

The bootstrap weights $W_{ni}$'s are assumed to belong to the
class of exchangeable bootstrap weights introduced in \cite{pw93}.
Specifically, they satisfy:

\begin{enumerate}[W2.]
\item[W1.] The vector $W_{n}=(W_{n1},\ldots,W_{nn})'$ is
exchangeable for all $n=1,2,\ldots,$ that is, for any permutation
$\pi=(\pi_{1},\ldots,\pi_{n})$ of $(1,2,\ldots,n)$, the joint
distribution of $\pi(W_{n})=(W_{n\pi_{1}},\ldots,W_{n\pi_{n}})'$
is the same as that of $W_{n}$.

\item[W2.] $W_{ni}\geq0$ for all $n$, $i$ and
$\sum_{i=1}^{n}W_{ni}=n$ for all $n$.

\item[W3.] For some positive constant $C < \infty$,
$\lim\sup_{n\rightarrow\infty}\|W_{n1}\|_{2,1}\leq C$, where
$\|W_{n1}\|_{2,1}=\int_{0}^{\infty}\sqrt{P_W(W_{n1}\geq u)}\,du$.

\item[W4.]
$\lim_{\lambda\rightarrow\infty}\lim\sup_{n\rightarrow\infty}
\sup_{t\geq\lambda}t^{2}P_W(W_{n1}>t)=0$.

\item[W5.]
$(1/n)\sum_{i=1}^{n}(W_{ni}-1)^{2}\stackrel{P_W}{\longrightarrow}c^{2}>0$.
\end{enumerate}

The bootstrap weights corresponding to Efron's nonparametric
bootstrap satisfy W1--W5. Another important class of bootstrap
whose weights satisfy W1--W5 is the \textit{multiplier bootstrap} in
which\vspace*{-1pt} $W_{ni}=\omega_i/\bar{\omega}_n$ and
$(\omega_1,\ldots,\omega_n)$ are i.i.d. positive r.v.s with
$\|\omega_1\|_{2,1}<\infty$. By taking\vspace*{1pt}
$\omega_i\stackrel{\mathrm{i.i.d.}}{\sim}\operatorname{Exp}(1)$, we obtain the
\textit{Bayesian bootstrap} of \cite{r81}. The multiplier bootstrap is
often thought to be a smooth alternative to the nonparametric
bootstrap \cite{l93}. In general, conditions W3--W5 are easily
satisfied under some moment conditions on $W_{ni}$; see Lemma 3.1
of \cite{pw93}. The sampling schemes that satisfy
conditions W1--W5 include \textit{the double bootstrap}, \textit{the urn
bootstrap and the grouped or delete-h Jackknife} \cite{e82}; see
\cite{pw93}. The value of $c$ in W5 is independent of $n$ and
depends on the resampling method, for example, $c=1$ for the
nonparametric bootstrap and Bayesian bootstrap, and $c=\sqrt{2}$
for the double bootstrap.

There exist two sources of randomness for the bootstrapped
quantity, for example, $\widehat{\theta}^{\ast}$ and
$\widehat{\eta}^\ast$: one comes from the observed data; another
comes from the resampling done by the bootstrap, that is, randomness in
$W_{ni}$'s. Therefore, in order to rigorously state our
theoretical results for the bootstrap, we need to specify relevant
probability spaces and define the related stochastic orders.

We view $X_{i}$ as the $i$th coordinate projection from the
canonical probability space $(\mathcal{X}^{\infty},
\mathcal{A}^{\infty}, P_X^{\infty})$ onto the $i$th copy of
$\mathcal{X}$. For the joint randomness involved, the product
probability space is defined as
\[
(\mathcal{X}^{\infty},\mathcal{A}^{\infty},P_{X}^{\infty})
\times(\mathcal{W},\Omega,P_{W})=(\mathcal{X}^{\infty}\times
\mathcal{W},\mathcal{A}^{\infty}\times\Omega,P_{\XW}).
\]
In this
paper, we assume that the bootstrap weights $W_{ni}$'s are
independent of the data $X_i$'s, thus $P_{\XW}=P_X^{\infty}\times
P_W$. We write $P_X^{\infty}$ as $P_X$ for simplicity thereafter.
Define $E_{\XW}^o$ as the outer expectation w.r.t. $P_{\XW}$. The
notation $E_{W|X}^o$, $E_X^o$ and $E_W$ are defined similarly.

Given a real-valued function $\Delta_n$ defined on the above
product probability space, for example, $\widehat{\theta}^{\ast}$,
we say
that $\Delta_n$ is of an order $o_{P_W}^o(1)$ in
$P_{X}^{o}$-probability if for any $\varepsilon,\delta>0$,
%
%e23 ###
%
\begin{equation}\label{bio}
P_{X}^{o}\{P_{W|X}^o(|\Delta_n|>\varepsilon)>\delta\}\longrightarrow
0 \qquad\mbox{as } n\rightarrow\infty,
\end{equation}
and that $\Delta_n$ is of an order $O_{P_W}^o(1)$ in
$P_{X}^{o}$-probability if for any $\delta>0$, there exists a
$0<M<\infty$ such that
%
%e24 ###
%
\begin{equation}\label{inter9}
P_{X}^{o}\{P_{W|X}^{o}(|\Delta_n|\geq M)> \delta\}\longrightarrow
0 \qquad\mbox{as } n\rightarrow\infty.
\end{equation}

Given a function $\Gamma_n$ defined only on
$(\mathcal{X}^{\infty},\mathcal{A}^{\infty},P_X^{\infty})$, if it
is of an order $o_{P_X}^{o}(1)$ [$O_{P_X}^{o}(1)$], then it is
also of an order $o_{P_{\XW}}^{o}(1)$ [$O_{P_{\XW}}^{o}(1)$] based
on the following argument:
\begin{eqnarray*}
P_{\XW}^{o}(|\Gamma_n|>\varepsilon)&=&E^{o}_{\XW}1\{|\Gamma_n|>\varepsilon
\}=
E_{X}E_{W|X}1\{|\Gamma_n|>\varepsilon\}^o\\
&=&E_X1\{|\Gamma_n|>\varepsilon\}^o=P_X^{o}\{|\Gamma_n|>\varepsilon\},
\end{eqnarray*}
where the third equation holds since $\Gamma_n$ does not depend on
the bootstrap weight. More results on transition of various
stochastic orders are given in Lemma \ref{term} of the \hyperref[app]{Appendix}.
Such results are used repeatedly in proving our bootstrap
consistency theorem.

To establish the bootstrap consistency, we need some additional
conditions. The first condition is the measurability condition,
denoted as $M(P_X)$. We say a class of functions $\mathcal{F}\in
M(P_X)$ if $\mathcal{F}$ possesses enough measurability so that
$\mathbb{P}_{n}$ can be randomized, that is, we can replace
$(\delta_{X_{i}}-P_X)$ by $(W_{ni}-1)\delta_{X_{i}}$, and Fubini's
theorem can be used freely. The detailed description for $M(P_X)$
is spelled out in \cite{gz90} and also given in the \hyperref[app]{Appendix} of
this paper. Define
$\mathcal{T}=\{\widetilde{m}(\theta,\eta)\dvtx\|\theta-\theta_{0}\|+\|
\eta
-\eta_{0}\|\leq
R\}$ for some $R>0$. For the rest of the paper, we assume
$\mathcal{T}\in M(P_X)$.

The second class of conditions parallels conditions S1--S3 used for
obtaining asymptotic normality of $\widehat{\theta}$ and is only
slightly stronger. Thus, the bootstrap consistency for $\theta$ is
almost automatically guaranteed once $\widehat{\theta}$ is shown
to be asymptotically normal. Let $S_{n}(x)$ be the envelop
function of the class $\mathcal{S}_{n} = \mathcal{S}_n(\delta_n)$
defined in (\ref{sndfn}), that is,
\[
S_{n}(x)=\sup_{\|\eta-\eta_{0}\|\leq\delta_{n}} \biggl|
\frac{\widetilde{m}(\theta_{0},\eta)-\widetilde{m}_{0}}
{\|\eta-\eta_{0}\|} \biggr|.
\]

The next condition controls the tail of this envelop function.

\begin{enumerate}[SB1.]
\item[SB1.] \textit{Tail probability condition}:
%
%e25 ###
%
\begin{equation}\label{addass}
\lim_{\lambda\rightarrow
\infty}\lim\sup_{n\rightarrow\infty}\sup_{t\geq\lambda}t^{2}P_X^{o}
\bigl(S_{n}(X_{1})>t\bigr)=0
\end{equation}
for any sequence $\delta_n\rightarrow0$.
\end{enumerate}

Let
$\dot{\mathcal{T}}=\{\partial\widetilde{m}(\theta,\eta)/\partial
\theta\dvtx
(\theta,\eta)\in\mathcal{C}_{n}\}$, where $\mathcal{C}_n =
\mathcal{C}_n (\delta_n)$ is defined in (\ref{eq:par-neighbor}).

\begin{enumerate}[SB2.]
\item[SB2.] We assume that $\dot{\mathcal{T}}\in M(P_X)\cap
L_{2}(P_X)$ and that $\dot{\mathcal{T}}$ is $P$-Donsker.
\end{enumerate}

Condition SB2 ensures that the size of the function class
$\dot{\mathcal{T}}$ is reasonable so that the bootstrapped
empirical processes
$\mathbb{G}_n^{\ast}\equiv\sqrt{n}(\mathbb{P}_n^{\ast}-\mathbb{P}_n)$
indexed by $\dot{\mathcal{T}}$ has a limiting process conditional
on the observations; see Theorem 2.2 in \cite{pw93}.

For any fixed $\theta$, define
\[
\widehat{\eta}_{\theta}^{\ast}=\mathop{\arg\max}_{\eta\in
\mathcal{H}}\mathbb{P}_n^{\ast}m(\theta,\eta).
\]
The next
condition says that $\widehat{\eta}_{\theta}^{\ast}$ should be
close to $\eta_0$ if $\theta$ is close to $\theta_0$.

\begin{longlist}[SB3.]
\item[SB3.] \textit{Bootstrap convergence rate condition}: there
exists a $\gamma\in(1/4,1/2]$ such that
%
%e26 ###
%
\begin{equation}\label{conrate-b}
\|\widehat{\eta}^{\ast}_{\widetilde{\theta}}-\eta_{0}\|=
O_{P_{W}}^{o}(\|\widetilde{\theta}-\theta_{0}\|\vee
n^{-\gamma}) \qquad\mbox{in } P_{X}^{o}\mbox{-probability}
\end{equation}
for any
$\widetilde{\theta}\stackrel{P_{\XW}^o}{\longrightarrow}\theta_0$.
\end{longlist}

Verifications of conditions SB1--SB2 will be discussed in
Section \ref{veri-con}. Two general theorems are given in
Section \ref{secrate} to aid verification of condition SB3.

Now we are ready to present our main results.
Theorem \ref{asythm-b} below says that the bootstrap distribution
of $(\sqrt{n}/c)(\widehat{\theta}^\ast-\widehat{\theta})$,
conditional on the observations, asymptotically imitates the
unconditional distribution of
$\sqrt{n}(\widehat{\theta}-\theta_0)$. Let $P_{W|\mathcal{X}_n}$
denote the conditional distribution given the observed data
$\mathcal{X}_n$.
\begin{theorem}\label{asythm-b}
Suppose\vspace*{-1pt} that $\widehat{\theta}$ and $\widehat{\theta}^{\ast}$
satisfy (\ref{exactsol}) and (\ref{exactsol-b}), respectively.
Assume that
$\widehat{\theta}\stackrel{P_X}{\longrightarrow}\theta_0$ and
$\widehat{\theta}^{\ast}\stackrel{P_{W}^{o}}
{\longrightarrow}\theta_0$ in $P_X^o$-probability. In addition,
assume that conditions~\textup{I}, \textup{S1--S3, SB1--SB3} and \textup{W1--W5}
hold. We have that
%
%e27 ###
%
\begin{equation}\label{bconratep}
\|\widehat{\theta}^{\ast}-\theta_{0}\|=O_{P_{W}}^{o}(n^{-1/2})
\end{equation}
in $P_{X}^{o}$-probability. Furthermore,
%
%e28 ###
%
\begin{equation}\label{bcons}
\sqrt{n}(\widehat{\theta}^{\ast}-\widehat{\theta})=
-A^{-1}\mathbb{G}_{n}^{\ast}\widetilde{m}_{0}+ o_{P_{W}}^{o}(1)
\end{equation}
in $P_{X}^{o}$-probability. Consequently,
%
%e29 ###
%
\begin{equation}\label{proconv}
\sup_{x\in\mathbb{R}^d} \bigl|P_{W|\mathcal{X}_n}\bigl(\bigl(\sqrt
{n}/c\bigr)(\widehat
{\theta}^{\ast}-
\widehat{\theta})\leq x\bigr)-P\bigl(N(0,\Sigma)\leq
x\bigr) \bigr|=o_{P_X}^o(1),
\end{equation}
where ``$\leq$'' is taken componentwise, $c$ is given in \textup{W5} and
$\Sigma\equiv A^{-1}B (A^{-1})'$ with $A$ and $B$ given in
(\ref{a}) and (\ref{b}), respectively. Thus, we have
%
%e30 ###
%
\begin{equation}\label{probc1}
%&&(\sqrt{n}/c)(\widehat{\theta}^{\ast}-\widehat{\theta})\Longrightarrow
%N(0,\Sigma) \mbox{in} P_X^o-\mbox{Probability},\label{bconcor}\\
\sup_{x\in\mathbb{R}^d} \bigl|P_{W|\mathcal{X}_n}\bigl(\bigl(\sqrt
{n}/c\bigr)(\widehat
{\theta}^{\ast}-\widehat{\theta})\leq
x\bigr)-P_X\bigl(\sqrt{n}(\widehat{\theta}-\theta_0)\leq x\bigr) \bigr|
\stackrel{P_X^o}{\longrightarrow}0.
\end{equation}
\end{theorem}

The consistency assumption for $\widehat{\theta}^\ast$ can be
established by adapting the Argmax theorem, that is, Corollary 3.2.3
in \cite{vw96}. Briefly, we need two conditions for accomplishing
this. The first one is the ``well-separated'' condition
(\ref{sepcon}). The second one is
%
%e31 ###
%
\begin{equation}\label{concon}
{\sup_{(\theta,\eta)\in\Theta\times\mathcal{H}}}|\mathbb
{P}_n^{\ast}
m(\theta,\eta)-P_X m(\theta,\eta)|\stackrel{P_{\XW}^{o}}
{\longrightarrow}0.
\end{equation}
By the multiplier Glivenko--Cantelli theorem, that is, Lemma 3.6.16 in
\cite{vw96}, and (\ref{term0}) in the \hyperref[app]{Appendix}, we know that
(\ref{concon}) holds if $\{m(\theta,\eta)\dvtx\theta\in\Theta,
\eta\in\mathcal{H}\}$ is shown to be $P$-Glivenko--Cantelli.
\begin{remark}\label{rm4}
For any consistent
$\widehat{\Sigma}^{\ast}\stackrel{P_{\XW}^o}{\longrightarrow}
\Sigma$ and
$\widehat{\Sigma}\stackrel{P_X}{\longrightarrow}\Sigma$, we have
%
%e32 ###
%
\begin{eqnarray}\label{probconv3}
&& \sup_{x\in\mathbb{R}^d} \bigl|P_{W|\mathcal{X}_n} \bigl(
\bigl(\sqrt{n}/c\bigr)(\widehat{\Sigma}^{\ast})^{-1/2}(\widehat{\theta
}^{\ast}
-\widehat{\theta})\leq x \bigr)\nonumber\\[-8pt]\\[-8pt]
&&\hspace*{56.15pt}{}
- P_X \bigl(\sqrt{n}\widehat{\Sigma}^{-1/2}(\widehat{\theta
}-\theta
_0)\leq
x \bigr) \bigr|\stackrel{P_X^o}{\longrightarrow}0\nonumber
\end{eqnarray}
by the arguments in proving Theorem \ref{asythm-b}, Slutsky's
theorem and Lemma \ref{term}. A possible candidate for the
consistent $\widehat{\Sigma}^{\ast}$ is the block jackknife
proposed in \cite{mk05b}.
\end{remark}
\begin{remark}
Our arguments in proving Theorem \ref{asythm-b} can also be used
to improve the remainder term in (\ref{bcons}) from
``$o_{P_W}^o(1)$ in $P_X^o$-probability'' to
``$O_{P_W}^o(n^{-2\gamma+1/2})$ in $P_X^o$-probability'' if we
strengthen the ``nearly maximizing'' condition (\ref{exactsol-b})
to the exactly maximizing condition (\ref{exactsol-b0}). A similar
result holds in Proposition \ref{asythm} where the remainder term
$o_{P_X}^{o} (1)$ in (\ref{asydis-nb}) can be improved to
$O_{P_X}(n^{-2\gamma+1/2})$ if (\ref{exactsol}) is strengthened to
(\ref{exactsol0}). It is interesting to note that the rate of
convergence of the remainder term depends on how accurately the
nuisance function parameter $\eta$ can be estimated. In
particular, if $\eta$ is $\sqrt{n}$-estimable, then the remainder
is of the order of $O(n^{-1/2})$.
\end{remark}

The distribution consistency result of the bootstrap estimator
$\widehat{\theta}^{\ast}$ proven in (\ref{probc1}) can be used to
prove the consistency of a variety of bootstrap confidence sets,
that is, \textit{percentile}, \textit{hybrid} and $t$ types.

A lower $\alpha$th quantile of bootstrap distribution is any
quantity $\tau_{n\alpha}^{\ast}\in\mathbb{R}^d$ satisfying
$\tau_{n\alpha}^{\ast}=\inf\{\varepsilon\dvtx
P_{W|\mathcal{X}_n}(\widehat{\theta}^{\ast}\leq\varepsilon)\geq
\alpha\}$,
where $\varepsilon$ is an infimum over the given set only if there
does not exist a $\varepsilon_1<\varepsilon$ in $\mathbb{R}^d$ such that
$P_{W|\mathcal{X}_n}(\widehat{\theta}^{\ast}\leq\varepsilon_1)\geq
\alpha$.
Because of the assumed smoothness of the criterion function
$m(\theta,\eta)$ in our setting, we can, without loss of
generality, assume $P_{W|\mathcal{X}_n}(\widehat{\theta}^{\ast}
\leq\tau_{n\alpha}^{\ast})=\alpha$. Due to the distribution
consistency result proven in (\ref{probc1}), we can approximate
the $\alpha$th quantile of the distribution of $(\widehat\theta-
\theta_0)$ by $(\tau_{n\alpha}^\ast-\widehat{\theta})/c$.
Thus, we define the \textit{percentile}-type bootstrap confidence set as
\[
\mathrm{BC}_{p}(\alpha)= \biggl[\widehat{\theta}+\frac{\tau_{n(
\alpha/2)}^\ast- \widehat{\theta}}{c},
\widehat{\theta}+\frac{\tau_{n(1-\alpha/2)}^\ast-
\widehat{\theta}}{c} \biggr].
\]
Similarly, we can approximate the $\alpha$th quantile of
$\sqrt{n}(\widehat{\theta}-\theta_0)$ by $\kappa_{n\alpha}^\ast$,
where $\kappa_{n \alpha}^{\ast}$ is the $\alpha$th quantile of the
hybrid quantity $(\sqrt{n}/c) (\widehat{\theta}^\ast-\widehat
{\theta})$,
that is,
$P_{W|\mathcal{X}_n}((\sqrt{n}/c)\times(\widehat{\theta}^{\ast}
-\widehat{\theta})\leq\kappa_{n \alpha}^{\ast})=\alpha$.
Thus, we define the \textit{hybrid}-type bootstrap confidence set as
\[
\mathrm{BC}_{h}(\alpha)= \biggl[\widehat{\theta}-
\frac{\kappa_{n(1-\alpha/2)}^\ast}{\sqrt{n}},\widehat{\theta}-
\frac{\kappa_{n(\alpha/2)}^\ast}{\sqrt{n}} \biggr].
\]
Note that
$\tau_{n\alpha}^{\ast}$ and $\kappa_{n\alpha}^{\ast}$ are not
unique since $\theta$ is assumed to be a vector.

We now prove the consistency of the above bootstrap confidence
sets by using the arguments in Lemma 23.3 of \cite{v98}. First, it
follows from (\ref{asynor}) and (\ref{proconv}) that, for any
$x\in\mathbb{R}^d$,
%
%e34 ###
%e33 ###
%
\begin{eqnarray}
\label{ab}
P_X\bigl(\sqrt{n}(\widehat{\theta}-\theta_0)\leq x\bigr)&\longrightarrow&
\Psi(x),\\
\label{ab2}
P_{W|\mathcal{X}_n}\bigl(\bigl(\sqrt{n}/c\bigr)(\widehat{\theta}^{\ast}-
\widehat{\theta})\leq x\bigr)&\stackrel{P_X^o}
{\longrightarrow}&\Psi(x),
\end{eqnarray}
where $\Psi(x)=P(N(0,\Sigma)\leq x)$. The quantile convergence
theorem, that is, Lemma 21.1 in \cite{v98}, applied to (\ref{ab2})
implies that $\kappa_{n\alpha}^{\ast}\to\Psi^{-1}(\alpha)$ almost
surely. When applying quantile convergence theorem, we use
the almost sure representation Theorem 2.19 in \cite{v98}
and argue along subsequences. Then the Slutsky's lemma
implies that $\sqrt{n}(\widehat{\theta}-\theta_0)-\kappa_{n(\alpha
/2)}^{\ast}$
weakly converges to $N(0,\Sigma)-\Psi^{-1}(\alpha/2)$. Thus,
\begin{eqnarray*}
P_{\XW} \biggl(\theta_0\leq\widehat{\theta}-
\frac{\kappa_{n(\alpha/2)}^{\ast}}{\sqrt{n}} \biggr)&=&P_{\XW} \bigl(
\sqrt{n}(\widehat{\theta}-\theta_0)\geq\kappa_{n(\alpha
/2)}^{\ast}
\bigr)\\
&\to&P_{\XW} \bigl(N(0,\Sigma)\geq
\Psi^{-1}(\alpha/2) \bigr)\\
&=& 1-\alpha/2.
\end{eqnarray*}
This argument yields the consistency of the \textit{hybrid}-type
bootstrap confidence set, that is, (\ref{perci1}) below, and can also be
applied to justify the \textit{percentile}-type bootstrap confidence
set, that is, (\ref{perci2}) below. The following Corollary \ref{perci}
summarizes the above discussion.
\begin{corollary}\label{perci}
Under the conditions in Theorem \ref{asythm-b}, we have
%
%e36 ###
%e35 ###
%
\begin{eqnarray}
\label{perci2}
P_{\XW} \bigl(\theta_0\in \mathrm{BC}_{p}(\alpha)
\bigr)&\longrightarrow&1-\alpha,\\
\label{perci1}
P_{\XW} \bigl(\theta_0\in
\mathrm{BC}_{h}(\alpha) \bigr)&\longrightarrow&1-\alpha
\end{eqnarray}
as $n\rightarrow\infty$.
\end{corollary}

It is well known that the above bootstrap confidence sets can be
computed easily through routine bootstrap sampling.

Investigating the consistency of the bootstrap variance estimator
is also of great interest. However, the usual sufficient condition
for moment consistency, that is, uniform integrability condition,
becomes very hard to verify due to the existence of an
infinite-dimensional parameter $\eta$. An alternative resampling
method to obtain the variance estimator in semiparametric models
is the block jackknife approach, which was proposed and
theoretically justified in \cite{mk05b}. We do not pursue this
topic further in this paper.
\begin{remark}
Provided consistent variance estimators $\widehat{\Sigma}^{\ast}$
and $\widehat{\Sigma}$ are available, we can define the
$t$-type bootstrap confidence set as
\[
\mathrm{BC}_{t}(\alpha)= \biggl[\widehat{\theta}-\frac{\widehat{\Sigma}^{1/2}
\omega_{n(1-\alpha/2)} ^{\ast}}{\sqrt{n}},\widehat{\theta}-\frac{
\widehat{\Sigma}^{1/2}
\omega_{n(\alpha/2)}^{\ast}}{\sqrt{n}} \biggr],
\]
where $\omega
_{n\alpha}^{\ast}$ satisfies
$P_{W|\mathcal{X}_n}((\sqrt{n}/c)(\widehat{\Sigma}^{\ast
})^{-1/2}(\widehat{\theta}
^{\ast}-\widehat{\theta})\leq\omega_{n\alpha}^{\ast})=\alpha$. By
applying again the arguments in Lemma 23.3 of \cite{v98} to
(\ref{asynorest}) and (\ref{probconv3}), we can prove that
\[
P_{\XW}\bigl(\theta_0\in \mathrm{BC}_{t}(\alpha)\bigr)\longrightarrow
1-\alpha
\]
as $n\rightarrow\infty$.
\end{remark}

%s4 ###
\section{Verifications of conditions S1, S2 and SB1, SB2}\label{veri-con}

%s4.1 ###
\subsection{Verifications of conditions \textup{S1} and \textup{S2}}\label{s2s3}
The continuity modulus condition (\ref{conmod}) in S1 can be
checked via one of the following two approaches. The first
approach is to show the boundedness of
$E^{o}_X\|\mathbb{G}_{n}\|_{\mathcal{S}_{n}}$ by using Lemma 3.4.2
in~\cite{vw96}. The second approach is to calculate the bracketing
entropy number of $\mathcal{S}_{n}$ and apply Lemma 5.13 in
\cite{v00} if $L_2$-norm is used on the nuisance parameter. As for
(\ref{conmod2}), we can verify it easily if we can show that the
class of functions
$\{(\partial/\partial\theta)\widetilde{m}(\theta,\eta)\dvtx
(\theta,\eta)\in\mathcal{C}_{n}\}$ is $P$-Donsker.

Next, we discuss how to verify the smoothness condition S2. We
first write $P_X(\widetilde{m}(\theta,\eta)-\widetilde{m}_0)$ as
the sum of
$P_X(\widetilde{m}(\theta,\eta)-\widetilde{m}(\theta_{0},\eta))$
and $P_X(\widetilde{m}(\theta_{0},\eta)-\widetilde{m}_{0})$. We
apply the Taylor expansion
to obtain
\begin{eqnarray*}
&&P_X\bigl(\widetilde{m}(\theta,\eta)-\widetilde{m}(\theta_{0},\eta)\bigr)\\
&&\qquad=P_X \{m_{11}(\theta_0,\eta)-m_{21}(\theta_0,\eta)[H^{\dag}
(\theta_0,\eta)] \}(\theta-\theta_{0})+O(\|\theta-\theta
_{0}\|
^{2})\\
&&\qquad=A(\theta-\theta_{0})+ (\theta-\theta_{0})O(\|\eta-\eta_{0}\|) +
O(\|\theta-\theta_0\|^2),
\end{eqnarray*}
where $A$ is defined in (\ref{a}), the first and second equality
follows from the Taylor expansion of $\theta\mapsto
P_X\widetilde{m}(\theta,\eta)$ around $\theta_0$ and
\[
\eta\mapsto
P_X \{m_{11}(\theta_0,\eta)-m_{21}(\theta_0,\eta)[H^{\dag}
(\theta_0,\eta)] \}
\]
around $\eta_0$, respectively. By
applying the second-order Taylor expansion to $\eta\mapsto
P_X\widetilde{m}(\theta_0,\eta)$ around $\eta_0$ and considering
(\ref{lfseq}), we can show that
$P(\widetilde{m}(\theta_{0},\eta)-\widetilde{m}_{0})=O(\|\eta-
\eta_0\|^2)$.
In summary, condition S2 usually holds in models where the map
$\eta\mapsto\widetilde{m}(\theta_0,\eta)$ is smooth in the sense
that the Fr\'{e}chet derivative of $\eta\mapsto
P_X((\partial/\partial\theta)\widetilde{m}(\theta_0,\eta))$ around
$\eta_0$ and the second order Fr\'{e}chet derivative of
$\eta\mapsto P_X\widetilde{m}(\theta_0,\eta)$ around $\eta_0$ are
bounded as discussed above.

%
%s4.2 ###
\subsection{Verifications of conditions \textup{SB1} and \textup{SB2}}\label{sb12}
We can verify\break condition SB1 by showing either $S_n(x)$ is
uniformly bounded, that is,\break $\lim\sup_{n\rightarrow\infty}
S_n(x)\leq
M<\infty$ for every $x\in\mathcal{X}$, or more generally,\break
$\lim\sup_{n\rightarrow\infty} E[\{S_n(X_1)\}^{2+\delta}]<\infty$
for some $\delta>0$. That the moment condition implies
condition SB1 follows from the Chebyshev's inequality. In our
examples in Section \ref{example}, the uniformly boundedness
condition is usually satisfied. Hence, we focus on how to show
$S_n(x)$ is uniformly bounded here. By the Taylor expansion in a
Banach space, we can write
$\widetilde{m}(\theta_0,\eta)-\widetilde{m}_0=D_{\widetilde{\eta}}
[\eta-\eta_0]$, where $\widetilde{\eta}$ lies on the line segment
between $\eta$ and $\eta_0$, and $D_{\xi}[h]$ is the Fr\'{e}chet
derivative of $\eta\mapsto\widetilde{m}(\theta_0,\eta)$ at $\xi$
along the direction $h$. Since we require
$\|\eta-\eta_0\|\leq\delta_n\rightarrow0$, the bounded
Fr\'{e}chet derivative at $\eta_0$ will imply that $S_n(x)$ is
uniformly bounded. The method in verifying (\ref{conmod2}) of
condition S1 can be applied to check condition SB2; see the
discussion in the previous subsection.

%
%s5 ###
\section{Convergence rates of bootstrap estimate of functional
parameter}\label{secrate} In this section, we present two general
theorems for calculating the convergence rate of the bootstrap
estimate of the functional parameter. These results can be applied
to verify condition SB3. Condition S3 can also be verified
based on these theorems by assuming the weights $W_{ni}=1$. Note
that both theorems extend general results on $M$-estimators
\cite{vw96,mv99} to bootstrap $M$-estimators and are also of
independent interest. Separate treatments are given to the cases
that the estimate $\eta$ has $\sqrt{n}$ convergence rate, that is,
Section \ref{rootn}, and has slower than $\sqrt{n}$ rate, that is,
Section~\ref{slow}.

%s5.1 ###
\subsection{Root-$n$ rate}\label{rootn}
We consider a collection of measurable
objective functions $x\mapsto k(\theta,\eta)[g](x)$ indexed by the
parameter $(\theta,\eta)\in\Theta\times\mathcal{H}$ and an
arbitrary index set $g\in\mathbf{G}$. For example,
$k(\theta,\eta)[g]$ can be the score function for $\eta$ given any
fixed $\theta$ indexed by $g\in\mathbf{G}$. Define
\begin{eqnarray*}
U_{n}^{\ast}(\theta,\eta)[g]&=&\mathbb{P}_{n}^{\ast}k(\theta,\eta
)[g],\\
U_{n}(\theta,\eta)[g]&=&\mathbb{P}_{n}k(\theta,\eta)[g],\\
U(\theta,\eta)[g]&=&P_Xk(\theta,\eta)[g].
\end{eqnarray*}
We assume that the maps $g\mapsto U_{n}^{\ast}(\theta,\eta)[g]$,
$g\mapsto U_{n}(\theta,\eta)[g]$ and $g\mapsto U(\theta,\eta)[g]$
are uniformly bounded, so that $U_{n}^{\ast}$, $U_{n}$ and $U$ are
viewed as maps from the parameter set $\Theta\times\mathcal{H}$
into $\ell^{\infty}(\mathbf{G})$. The following conditions are
assumed in Theorem~\ref{ratethm2} below:
%
%e38 ###
%e37 ###
%
\begin{equation}
\label{ratecon1}
\{k(\theta,\eta)[g]\dvtx
\|\theta-\theta_{0}\|+\|\eta-\eta_{0}\|\leq\delta,
g\in\mathbf{G}\}\in M(P_X)\cap
L_2(P_X)
\end{equation}
and is $P$-Donsker for some $\delta>0$,
\begin{equation}
\label{ratecon2}
\sup_{g\in\mathbf{G}}P_X \{k(\theta,\eta)[g]-k(\theta
_{0},\eta_{0}
)[g] \}^{2}\rightarrow
0 \qquad\mbox{as } \|\theta-\theta_{0}\|+
\|\eta-\eta_0\|\rightarrow0.\hspace*{-28pt}
\end{equation}
Let
\[
\mathcal{D}_{n}= \biggl\{\frac{k(\theta,\eta)[g]
-k(\theta_0,\eta_0)[g]}{1+\sqrt{n}\|\theta-\theta_0\|+\sqrt{n}\|
\eta
-\eta_0\|}\dvtx g\in
\mathbf{G},\|\theta-\theta_0\|+\|\eta-\eta_0\|\leq\delta
_{n} \biggr\}
\]
and $D_{n}(X)$ be the envelop function of the class of functions
$\mathcal{D}_n$. For any sequence $\delta_n\rightarrow0$, we
assume that $D_n(X)$ satisfies
%
%e39 ###
%
\begin{equation}\label{ratecon7}
\lim_{\lambda\rightarrow\infty}\lim\sup_{n\rightarrow\infty
}\sup_{t
\geq\lambda}t^{2}P_{X}^{o}\bigl(D_{n}(X_1)>t\bigr)=0.
\end{equation}

Now we consider the convergence rate of
$\widehat{\eta}_{\widetilde{\theta}}^{\ast}$ satisfying:
%
%e40 ###
%
\begin{equation}\label{neunb-b}
U_{n}^{\ast}(\widetilde{\theta},\widehat{\eta}_{\widetilde{\theta
}}^{\ast})[g]=O_{P_{\XW}}^{o}(n^{-1/2})
\end{equation}
for any
$\widetilde{\theta}\stackrel{P_{\XW}^o}{\longrightarrow}\theta_0$
and $g$ ranging over $\mathbf{G}$. In Theorem \ref{ratethm2}
below, we will show that
$\widehat{\eta}_{\widetilde{\theta}}^{\ast}$ has the root-$n$
convergence rate under conditions (\ref{ratecon1})--(\ref{ratecon7}).
\begin{theorem}\label{ratethm2}
Suppose that $U\dvtx
\Theta\times\mathcal{H}\mapsto\ell^{\infty}(\mathbf{G})$ is
Fr\'{e}chet differentiable at $(\theta_{0},\eta_{0})$ with bounded
derivative $\dot{U}\dvtx\mathbb{R}^{d}\times
\operatorname{lin}\mathcal{H}\mapsto\ell^{\infty}(\mathbf{G})$ such that
the map
$\dot{U}(0,\cdot)\dvtx\operatorname{lin}\mathcal{H}\mapsto\ell^{\infty}
(\mathbf{G})$ is invertible with an inverse that is continuous on
its range. Furthermore, assume that
(\ref{ratecon1})--(\ref{ratecon7}) hold, and that
$U(\theta_{0},\eta_{0})=0$, then
%
%e41 ###
%
\begin{equation}\label{nuirate}
\|\widehat{\eta}^{\ast}_{\widetilde{\theta}}-\eta_{0}\|
=O_{P_{W}}^{o}(\|
\widetilde{\theta}-\theta_{0}\|\vee
n^{-1/2})
\end{equation}
in $P_X^o$-probability, given that
$\widetilde{\theta}\stackrel{P_{\XW}^o}{\longrightarrow}\theta_0$
and $\widehat{\eta}^{\ast}_{\widetilde{\theta}}
\stackrel{P_{\XW}^o}{\longrightarrow}\eta_0$.
\end{theorem}

The proof of Theorem \ref{ratethm2} is given in Appendix \ref{appA4}.

%s5.2 ###
\subsection{Slower than root-$n$ rate}\label{slow}
We next present a result that deals with slower than
$\sqrt{n}$ convergence rate for the bootstrap $M$-estimate of the
functional parameter. This result is so general that it can be
applied to the sieve estimate of nuisance parameter
\cite{g81}. The essence of the sieve method is that a sequence of
increasing spaces (sieves), that is, $\mathcal{H}_n$, is employed to
approximate the large parameter space, for example, $\mathcal{H}$. In
other words, for any $\eta\in\mathcal{H}$, there exists a
$\pi_{n}\eta\in\mathcal{H}_n$ such that
$\|\eta-\pi_{n}\eta\|\rightarrow0$ as $n\rightarrow\infty$.

Now, we consider the $M$-estimate
$\widehat{\eta}_{\theta}^{\ast}\in\mathcal{H}_{n}$ satisfying
%
%e42 ###
%
\begin{equation}\label{mconditi}
\mathbb{P}_{n}^{\ast}v(\theta,\widehat{\eta}^{\ast}_{\theta
})\geq
\mathbb{P}_{n}^{\ast}v(\theta,\eta_{n}) \qquad\mbox{for
any } \theta\in\Theta\mbox{ and some
$\eta_{n}\in\mathcal{H}_{n}$},
\end{equation}
where $x\mapsto v(\theta,\eta)(x)$ is a measurable objective
function. Let ``$\lesssim$'' and ``$\gtrsim$'' denote greater than or
smaller than, up to an universal constant. We assume the following
conditions hold for every $\delta>0$:
%
%e43 ###
%
\begin{eqnarray}\label{ratecon3}
&\displaystyle E_X\bigl(v(\theta,\eta)-v(\theta,\eta_{n})\bigr)\lesssim-d^{2}
(\eta,\eta_{n})+\|\theta-\theta_{0}\|^{2},&
\\
%
%e44 ###
%
\label{ratecon4}
&\hspace*{4.5pt}\displaystyle E_X^{o}\sup_{\theta\in\Theta,\eta\in\mathcal{H}_{n},
\|\theta-\theta_{0}\|\leq\delta,d(\eta,\eta_{n})\leq\delta}
\bigl|\mathbb{G}_{n}\bigl(v(\theta,\eta)-v(\theta,\eta_{n})\bigr)\bigr|
\lesssim\psi_{n}(\delta),&
\\
%
%e45 ###
%
\label{ratecon5}
&\displaystyle {E_{\XW}^{o}\sup_{\theta\in\Theta,\eta
\in\mathcal{H}_{n},\|\theta-\theta_{0}\|\leq\delta,d(\eta,\eta
_{n})\leq
\delta}}
\bigl|\mathbb{G}_{n}^{\ast}\bigl(v(\theta,\eta)-v(\theta, \eta
_{n})\bigr)\bigr|\lesssim
\psi_{n}^{\ast}(\delta).&
\end{eqnarray}
Here $d^2(\eta,\eta_n)$ may be thought of as the square of a
distance, for example, $\|\eta-\eta_n\|^2$, but our theorem is also true
for any arbitrary function $\eta\mapsto d^2(\eta,\eta_n)$.
%In the above, we assume the distance function
%$d(\cdot,\cdot)$ may depend on $\theta$. This allows us
%to compute the convergence rate for the penalized nuisance
%parameter in which the smoothing parameter can be included in
%$\theta$.
%
\begin{theorem}\label{ratethm}
Suppose that conditions (\ref{ratecon3})--(\ref{ratecon5}) hold. We
assume (\ref{ratecon4}) [and (\ref{ratecon5})] is valid for\vspace*{1pt} functions
$\psi_{n}$ $(\psi_{n}^{\ast})$ such that
$\delta\mapsto\psi_{n}(\delta)/\delta^{\alpha}$
[$\delta\mapsto\psi_{n}^{\ast}(\delta)/\delta^{\alpha}$] is
decreasing for some $0<\alpha<2$. Then for every
$(\widetilde{\theta},\widehat{\eta}_{\widetilde{\theta}}^{\ast})$
satisfying
$P(\widetilde{\theta}\in\Theta,\widehat{\eta}_{\widetilde{\theta}}^
{\ast}\in\mathcal{H}_{n})\rightarrow1$, we have
\[
d(\widehat{\eta}_{\widetilde{\theta}}^{\ast},\eta_{n})\leq
O_{P_{W}}^{o}(\delta_{n}\vee\|\widetilde{\theta}-\theta_{0}\|)
\]
in $P_X^o$-probability, for any sequence of positive numbers
$\delta_{n}$ satisfying both
$\psi_{n}(\delta_{n})\leq\sqrt{n}\delta_{n}^{2}$ and
$\psi_{n}^{\ast}(\delta_{n})\leq\sqrt{n}\delta_{n}^{2}$ for large
$n$.
\end{theorem}

The proof of Theorem \ref{ratethm} is given in Appendix \ref{appA5}.

In application of Theorem \ref{ratethm}, the parameter $\eta_n$ is
taken to be some element in $\mathcal{H}_n$ that is very close to
$\eta_0$. When $\mathcal{H}_{n}=\mathcal{H}$, a natural choice for
$\eta_{n}$ is $\eta_{0}$ and we can directly use
Theorem \ref{ratethm} to derive the convergence rate
$d(\widehat{\eta}_{\widetilde{\theta}}^{\ast},\eta_{0})$ as shown
in the examples of Section \ref{example}. In general, $\eta_{n}$
may be taken as the maximizer of the mapping $\eta\mapsto P_X
v(\theta_{0},\eta)$ over $\mathcal{H}_{n}$, the projection of
$\eta_{0}$ onto $\mathcal{H}_{n}$. Then we need to consider the
approximation rate of the sieve space $\mathcal{H}_{n}$ to
$\mathcal{H}$, that is, $d(\eta_{n},\eta_{0})$, since
$d(\widehat{\eta}_{\widetilde{\theta}}^{\ast},\eta_{0})\leq
d(\widehat{\eta}_{\widetilde{\theta}}^{\ast},\eta_{n})+
d(\eta_{n},\eta_{0})$. The approximation rate
$d(\eta_{n},\eta_{0})$ depends on the choices of sieves and is
usually derived in the mathematical literature. %Here, we skip
%further discussion and refer readers to \cite{c07}.

Now we discuss verification of the nontrivial conditions
(\ref{ratecon3})--(\ref{ratecon5}). The smoothness condition for
$v(\theta,\eta)$, that is, (\ref{ratecon3}), is implied by
%
%e47 ###
%e46 ###
%
\begin{eqnarray}
\label{smo1}
E_X\bigl(v(\theta,\eta)-v(\theta_{0},\eta_{n})\bigr)
&\lesssim& -d^{2}(\eta,\eta_{n})-\|\theta-\theta_{0}\|^{2},\\
\label{smo2}
E_X\bigl(v(\theta,\eta_{n})-v(\theta_{0},\eta_{n})\bigr)&\gtrsim&-\|\theta
-\theta
_{0}\|^{2}.
\end{eqnarray}
The two conditions depict the quadratic behaviors of the criterion
functions $(\theta,\eta)\mapsto E_Xv(\theta,\eta)$ and
$\theta\mapsto E_Xv(\theta,\eta_{n})$ around the maximum point
$(\theta_{0},\eta_{n})$ and~$\theta_{0}$, respectively. We next
present one useful lemma for verifying the continuity modulus of
(bootstrapped) empirical processes, that is,
(\ref{ratecon4}) and (\ref{ratecon5}). Denote
%
%e48 ###
%
\begin{equation}\label{addass2}
\mathcal{V}_{\delta}=\{x\mapsto
[v(\theta,\eta)(x)-v(\theta,\eta_{n})(x)]\dvtx
d(\eta,\eta_{n})\leq\delta,
\|\theta-\theta_{0}\|\leq\delta\}
\end{equation}
and define the bracketing entropy integral of
$\mathcal{V}_{\delta}$ as
%
%e49 ###
%
\begin{equation}\label{thm2k}
K(\delta,\mathcal{V}_{\delta},L_{2}(P_X))=\int_{0}^{\delta}\sqrt
{1+\log
N_{[\cdot]}(\varepsilon,\mathcal{V}_{\delta},L_{2}(P_X))}\,d\varepsilon,
\end{equation}
where $\log N_{[\cdot]}(\delta, \mathcal{A},d)$ is the
$\delta$-bracketing entropy number for the class $\mathcal{A}$
under the distance measure $d$.
\begin{lemma}\label{lerate}
Suppose that the functions $(x,\theta,\eta)\mapsto
v_{\theta,\eta}(x)$ are uniformly bounded for $(\theta,\eta)$
ranging over some neighborhood of $(\theta_{0},\eta_{n})$ and that
%
%e50 ###
%
\begin{equation}\label{lerate1}
E_X(v_{\theta,\eta}-v_{\theta,\eta_{n}})^{2}\lesssim
d^{2}(\eta,\eta_{n})+\|\theta-\theta_{0}\|^{2}.
\end{equation}
Then condition (\ref{ratecon4}) is satisfied for any functions
$\psi_{n}$ such that
%
%e51 ###
%
\begin{equation} \label{verlem1}
\psi_{n}(\delta)\geq K(\delta,\mathcal{V}_{\delta},L_{2}(P_X))
\biggl(1+\frac{K(\delta,\mathcal{V}_{\delta},L_{2}(P_X))}
{\delta^{2}\sqrt{n}} \biggr).
\end{equation}
Let $V_{n}(X)$ be the envelop function of the class
$\mathcal{V}_{\delta_n}$. If we further assume that, for each
sequence $\delta_{n}\rightarrow0$, the envelop functions $V_{n}$
satisfies
%
%e52 ###
%
\begin{equation}\label{inter15}
\lim_{\lambda\rightarrow\infty}\lim\sup_{n\rightarrow\infty}
\sup_{t\geq\lambda}t^{2}P_{X}^{o}\bigl(V_{n}(X_{1})>t\bigr)=0,
\end{equation}
then condition (\ref{ratecon5}) is satisfied for any functions
$\psi_{n}^{\ast}$ such that
%
%e53 ###
%
\begin{equation}\label{verlem2}
\psi_{n}^{\ast}(\delta)\geq
K(\delta,\mathcal{V}_{\delta},L_{2}(P_X))
\biggl(1+\frac{K(\delta,\mathcal{V}_{\delta},L_{2}(P_X))}
{\delta^{2}\sqrt{n}} \biggr).
\end{equation}
\end{lemma}
\begin{remark}
Note that the inequalities
$\psi_{n}(\delta)\lesssim\sqrt{n}\delta^{2}$ and
$\psi_{n}^{\ast}(\delta)\lesssim\sqrt{n}\delta^{2}$ are equivalent
to $K(\delta, \mathcal{V}_{\delta},
L_2(P_X))\lesssim\sqrt{n}\delta^{2}$ when we let $\psi_{n}$ and
$\psi_{n}^{\ast}$ be equal to the right-hand side of (\ref{verlem1})
and (\ref{verlem2}), respectively. Consequently, the convergence
rate of $\widehat{\eta}_{\widetilde{\theta}}^{\ast}$ calculated in
Theorem \ref{ratethm}, that is, $\delta_{n}$, is determined by the
bracketing entropy integral of $\mathcal{V}_{\delta_n}$.
\end{remark}
\begin{remark}
The assumptions of Lemma \ref{lerate} are relaxable to great
extent. For example, we can drop the uniform bounded condition on
the class of functions $v(\theta,\eta)$ by using the ``Bernstein
norm,'' that is, $\|f\|_{P,B}=(2P(e^{|f|}-1-|f|))^{1/2}$, instead of
the $L_{2}$-norm. In some cases, the bracketing entropy integral
diverges at zero. Then we can change the limit of the integration in
(\ref{thm2k}) from $[0,\delta]$ to $[a\delta^2\wedge\delta/3,
\delta]$
for some
small positive constant $a$, see Lemma 3.4.3 and page 326 in
\cite{vw96}.
\end{remark}

%s6 ###
\section{Examples}\label{example}
In this section, we apply the main results in
Section \ref{bootcons} to justify the bootstrap validity of
drawing semiparametric inferences in three examples of
semiparametric models. In the Cox regression models with censored
data, we use the log-likelihood as the criterion function, while
in the partially linear model, the least squares criterion is
used. The $M$-estimate of the nuisance functional parameters have
different convergence rates in these examples. Indeed, the
advantages of using bootstrap approach in all of the three
examples were considered in the literature, for example,
\cite{et86,lhs00}. This section also serves the purpose of illustration on
verification of the technical conditions used in the general
results.

%s6.1 ###
\subsection{Cox regression model with right censored data}
In the Cox regression model, the hazard function of the survival
time $T$ of a subject with covariate $Z$ is modeled as
%
%e54 ###
%
\begin{equation}\label{eg1den}
\lambda(t|z)\equiv\lim_{\Delta\rightarrow
0}\frac{1}{\Delta}P(t\leq T<t+\Delta|T\geq
t,Z=z)=\lambda(t)\exp(\theta' z),
\end{equation}
where $\lambda$ is an unspecified baseline hazard function and
$\theta$ is a regression vector. In this model, we are usually
interested in $\theta$ while treating the cumulative hazard
function $\eta(y)=\int_{0}^{y}\lambda(t)\,dt$ as the nuisance
parameter. The MLE for $\theta$ is proven to be semiparametric
efficient and widely used in applications. Here we consider
bootstrapping $\widehat{\theta}$, which corresponds to treating
log-likelihood as the criterion function $m(\theta,\eta)$ in our
general formulation.

With right censoring of survival time, the data observed is
$X=(Y,\delta,Z)$, where $Y=T\wedge C$, $C$ is a censoring time,
$\delta=I\{T\leq C\}$, and $Z$ is a regression covariate belonging
to a compact set $\mathbb{Z}\subset\mathbb{R}^d$. We assume that
$C$ is independent of $T$ given~$Z$. The log-likelihood is
obtained as
%
%e55 ###
%
\begin{equation}
m(\theta,\eta)=\delta\theta' z-\exp(\theta'
z)\eta(y)+\delta\log\eta\{y\},
\end{equation}
where $\eta\{y\} = \eta(y) - \eta(y-)$ is a point mass that
denotes the jump of $\eta$ at point~$y$. The parameter space
$\mathcal{H}$ is restricted to a set of nondecreasing cadlag
functions on the interval $[0,\tau]$ with $\eta(\tau)\leq M$ for
some constant $M$. By some algebra, we have\looseness=1
\begin{eqnarray*}
\widetilde{m}(\theta,\eta)(x)&=&m_{1}(\theta,\eta)-m_{2}(\theta
,\eta)
[H^{\dag}(\theta,\eta)]\\
&=& [\delta z-z\exp(\theta' z)\eta(y) ]\\
&&{} - \biggl[\delta
H^{\dag}(\theta,\eta)(y)-\exp(\theta'
z)\int_{0}^{y}H^{\dag}(\theta,\eta)(u)\,d\eta(u) \biggr],
\end{eqnarray*}
where
\[
H^{\dag}(\theta,\eta)(y)=\frac{E_{\theta,\eta}Z\exp(\theta'
Z)1\{Y\geq y\}}{E_{\theta,\eta}\exp(\theta' Z)1\{Y\geq y\}}.
\]

Conditions I, S1--S3 in guaranteeing the asymptotic normality of
$\widehat{\theta}$ have been verified in \cite{ck08a}. In
particular, the convergence rate of the estimated nuisance
parameter is established in Theorem 3.1 of \cite{mv99}, that is,
%
%e56 ###
%
\begin{equation}\label{eg1crate}
\|\widehat{\eta}_{\widetilde{\theta}_{n}}
-\eta_{0}\|_{\infty}=O_{P_X}(n^{-{1/2}}+
\|\widetilde{\theta}_{n}-\theta_{0}\|),
\end{equation}
where \mbox{$\|\cdot\|_{\infty}$} denotes the supreme norm. We next
verify the bootstrap consistency conditions, that is, SB1--SB3.
Condition SB1 trivially holds since it is easy to show that
$\eta\mapsto\widetilde{m}(\theta_0,\eta)$ has bounded Fr\'{e}chet
derivative around $\eta_0$. The $P$-Donsker condition SB2 has been
verified when verifying (\ref{conmod2}) in condition S1. In the
end, we will verify the bootstrap convergence rate condition
$\|\widehat{\eta}^{\ast}_{\widetilde{\theta}}-\eta_0\|_{\infty
}=O_{P_{\XW}}^{o}
(\|\widetilde{\theta}-\theta_0\|\vee n^{-1/2})$ via
Theorem \ref{ratethm2}. Since $\widehat{\eta}_{\theta}^{\ast}$
maximizes $\mathbb{P}_n^{\ast}m(\theta,\eta)$ for fixed~$\theta$,
we set $k(\theta,\eta)[g]=m_2(\theta,\eta)[g]$ and have
$U_n^{\ast}(\theta,\widehat{\eta}_{\theta}^{\ast})[g]=\mathbb
{P}_n^{\ast}
m_{2}(\theta,\widehat{\eta}_{\theta}^{\ast})[g]=0$. The
invertibility of $\dot{W}(0,\cdot)$, conditions
(\ref{ratecon1}) and (\ref{ratecon2}) have been verified in
\cite{mv99} when they showed (\ref{eg1crate}). Now we only need
to consider condition (\ref{ratecon7}): for $n$ so large that
$\delta_n\leq R$
\begin{eqnarray*}
D_n(x)&\equiv&\sup\biggl\{\frac{|(m_2(\theta,\eta)[g])-m_2(\theta_0,
\eta_0)[g]|}
{1+\sqrt{n}(\|\theta-\theta_0\|+\|\eta-\eta_0\|_{\infty})},
g\in\mathbf{G},\\
&&\hspace*{77pt}
\|\theta-\theta_0\|+\|\eta-
\eta_0\|_{\infty}\leq
\delta_n \biggr\}\\
&\leq&2\sup\{|m_2(\theta,\eta)[g]|,
g\in\mathbf{G},\|\theta-\theta_0\|+\|\eta-\eta_0\|_{\infty}\leq
R \}\\
&\leq&\mbox{some constant}.
\end{eqnarray*}
The last inequality follows from the assumption that $\mathbf{G}$
is a class of functions of bounded total variation and the
inequality that
$\int_{0}^{y}g(u)\,d\eta(u)\leq\eta(\tau)\|g\|_{\mathrm{BV}}$, where
$\|g\|_{\mathrm{BV}}$ is the total variation of the function $g$. Thus, condition (\ref{ratecon7}) holds trivially.

%s6.2 ###
\subsection{Cox regression model with current status data} We
next consider the current status data when each subject is
observed at a single examination time $C$ to determine if an event
has occurred. The event time $T$ cannot be known exactly. Then the
observed data are $n$ i.i.d. realizations of $X=(C, \delta, Z)\in
R^{+}\times\lbrace0,1 \rbrace\times\mathbb{Z}$, where
$\delta=I\{T \leq C\}$. The corresponding criterion function,
that is, the log-likelihood, is derived as
%
%e57 ###
%
\begin{equation}\label{eg1lik}
m(\theta,\eta) = \delta\log[1-\exp(-\eta(c)\exp(\theta'
z))]-(1-\delta)\exp(\theta' z)\eta(c).
\end{equation}
We make the following assumptions throughout the rest of this
subsection: (i)~$T$ and $C$ are independent given $Z$; (ii) the
covariance of $Z-E(Z|C)$ is positive definite, which guarantees
the efficient information to be positive definite; (iii)~$C$
possesses a Lebesgue density which is continuous and positive on
its support $[\sigma,\tau]$, for which the true nuisance parameter
$\eta_0$ satisfies $\eta_0(\sigma-)>0$ and
$\eta_0(\tau)<M<\infty$, and this density is continuously
differentiable on $[\sigma,\tau]$ with derivative bounded above
and bounded below by zero. The form of
$\widetilde{m}(\theta,\eta)$ can be found in \cite{ck08b} as
follows
\begin{eqnarray*}
\widetilde{m}(\theta,\eta)&=&m_{1}(\theta,\eta)-m_{2}(\theta,\eta)
[H^{\dag}(\theta,\eta)]\\
&=&\bigl(z\eta(c)-H^{\dag}(\theta,\eta)(c)\bigr)Q(x;\theta,\eta),
\end{eqnarray*}
where
\[
Q(x;\theta,\eta)=e^{\theta'z} \biggl[\frac{\delta}{\exp(e^{\theta
'z}\eta(c)
)-1}-(1-\delta) \biggr]
\]
and the form of
$H^{\dag}(\theta,\eta)(c)$ is given in (4) of \cite{ck08b}.

Conditions I and S1--S3 are verified in \cite{ck08b}. Conditions
SB1 and SB2 can be checked similarly as in the previous example.
Note that the convergence rate for the nuisance parameter becomes
slower, that is,
%
%e58 ###
%
\begin{equation}\label{eg1rate}
\|\widehat{\eta}_{\widetilde{\theta}_{n}}-\eta_{0}\|_{2}=
O_{P_X}(\|\widetilde{\theta}_{n}-\theta_{0}\|+n^{-1/3}),
\end{equation}
where \mbox{$\|\cdot\|_2$} denotes the regular $L_2$-norm, as shown in
\cite{mv99}. By Theorem \ref{ratethm}, we can show that the same
convergence rate, that is, $n^{-1/3}$, also holds for
$\widehat{\eta}_{\theta}^{\ast}$. The assumptions (\ref{ratecon3})
and (\ref{ratecon4}) in Theorem \ref{ratethm} are verified in
\cite{mv99} when showing (\ref{eg1rate}). We apply
Lemma \ref{lerate} to verify assumption (\ref{ratecon5}). We
show that condition (\ref{inter15}) on the envelop function
$V_n(x)$ holds: for $n$ so large that $\delta_n\leq R$
\begin{eqnarray*}
V_n(x)&\equiv&\sup\{|m(\theta,\eta)-m(\theta,\eta_0)|\dvtx
\|\eta-\eta_0\|_{2}\leq\delta_n, \|\theta-\theta_0\|\leq\delta
_n \}
\\
&\leq&2\sup\{|m(\theta,\eta)|\dvtx\|\eta-\eta_0\|_{2}\leq R,
\|\theta-\theta_0\|\leq R\}\\
&\leq&\mbox{some constant}.
\end{eqnarray*}

%s6.3 ###
\subsection{Partially linear models}\label{partl}
In this example, a continuous outcome variable~$Y$, depending on
the covariates $(W,Z)\in[0,1]^2$, is modeled as
\[
Y=\theta W+f(Z)+\xi,
\]
where $\xi$ is independent of $(W,Z)$ and
$f$ is an unknown smooth function belonging to
$\mathcal{H}\equiv\{f\dvtx[0,1]\mapsto[0,1],
\int_{0}^{1}(f^{(k)}(u))^2\,du\leq M\}$ for a fixed $0<M<\infty$. In
addition, we assume $E(\operatorname{Var}(W|Z))$ is positive definite and
$E\{f(Z)\}=0$. We want to estimate $(\theta,f)$ using the least
square criterion:
%
%e59 ###
%
\begin{equation}\label{plm}
m(\theta,f)=-\bigl(y-\theta w-f(z)\bigr)^2.
\end{equation}
Note that the above model would be more flexible if we did not
require knowledge of $M$. A sieve estimator could be obtained if
we replaced $M$ with a sequence $M_n\rightarrow\infty$. The theory
we develop in this paper will be applicable in this setting, but,
in order to maintain clarity of exposition, we have elected not to
pursue this more complicated situation here. Another approach is
to use penalization, the study of which is beyond the scope of
this paper.

Simple calculations give
\begin{eqnarray*}
\widetilde{m}(\theta,\eta)(x)&=&m_1(\theta,\eta)-m_2(\theta,\eta)
[H^{\dag}(\theta,\eta)]\\
&=&2\bigl(y-\theta w-f(z)\bigr)\bigl(w-H^{\dag}(\theta,\eta)(z)\bigr),
\end{eqnarray*}
where
\[
H^{\dag}(\theta,\eta)(z)=\frac{E_{\theta,\eta}(W(Y-\theta
W-f(Z))^2|Z=z)}{E_{\theta,\eta}((Y-\theta W-f(Z))^2|Z=z)}.
\]
The
finite variance condition I follows from
$E[W\{W-H^{\dag}(\theta_0,\eta_0)(Z)\}]>0$. The distribution of
$\xi$ is assumed to have finite second moment and satisfy
(\ref{zeroeq}), for example, $\xi\sim N(0,1)$. Conditions S1--S3 and SB2
can be verified using similar arguments in Example 3 of
\cite{ck08b}, in particular,
$\|\widehat{f}_{\widetilde{\theta}}-f_0\|_{2}=O_{P_X}(\|\widetilde
{\theta}-\theta_0\|\vee n^{-k/(2k+1)})$ in (\ref{conrate}). It is
easy to show that the Fr\'{e}chet derivative of
$\eta\mapsto\widetilde{m}(\theta_0,\eta)$ is bounded around
$\eta_0$, and thus the tail condition SB1 holds. To prove
$\|\widehat{f}_{\widetilde{\theta}}^{\ast}-f_0\|_2=O_{P_{\XW}}^{o}
(\|\widetilde{\theta}-\theta_0\|\vee n^{-k/(2k+1)})$ via
Theorem \ref{ratethm}, we proceed as in the previous example,
checking assumption (\ref{inter15}) using similar arguments,
that is, $V_n(x)$ is uniformly bounded.

%s7 ###
\section[Proof of Theorem 1 (bootstrap consistency
theorem)]{Proof of Theorem \protect\ref{asythm-b} (bootstrap consistency
theorem)}\label{pfthm2}

To prove Theorem~\ref{asythm-b}, we need the following lemma whose
proof is given in Appendix \ref{appA3}.
\begin{lemma}\label{bmodu}
Under the assumptions of Theorem \ref{asythm-b}, we have
%
%e60 ###
%
\begin{equation}\label{conmod-b}
\mathbb{G}_{n}^{\ast} \bigl(\widetilde{m}(\theta,\eta)-\widetilde
{m}(\theta_{0},\eta_{0}) \bigr)=
O_{P_{W}}^{o}(\|\theta-\theta_{0}\| \vee\|\eta-
\eta_{0}\|)
\end{equation}
in $P_{X}^{o}$-probability for $(\theta,\eta)\in\mathcal{C}_{n}$.
\end{lemma}

We shall use repeatedly Lemma \ref{term} in the \hyperref[app]{Appendix}, which
concerns about the transition of stochastic orders among different
probability spaces.

We first prove (\ref{bconratep}). Recall that
$\mathbb{G}_n=\sqrt{n}(\mathbb{P}_n-P_X)$ and
$\mathbb{G}_n^{\ast}=\sqrt{n}(\mathbb{P}_n^{\ast}-\mathbb{P}_n)$.
Define $\widehat{m}^{\ast}$ as
$\widetilde{m}(\widehat{\theta}^{\ast}, \widehat{\eta}^{\ast})$.
By some algebra, we have
\begin{eqnarray*}
&&\mathbb{G}_{n}^{\ast}\widetilde{m}_{0}+\mathbb{G}_{n}\widetilde
{m}_{0}+\sqrt{n}
P_X(\widehat{m}^{\ast}
-\widetilde{m}_{0})\\
&&\qquad=\mathbb{G}_{n}^{\ast}(\widetilde{m}_{0}-\widehat{m}^{\ast
})+\mathbb
{G}_{n}(\widetilde{m}_{0}-
\widehat{m}^{\ast})+\sqrt{n}\mathbb{P}
_{n}^{\ast}\widehat{m}^{\ast},
\end{eqnarray*}
since $P_X\widetilde{m}_{0}=0$. Thus, we have the following inequality:
%
%e61 ###
%
\begin{eqnarray}\label{equarat}
\bigl\|\sqrt{n}P_X(\widehat{m}^{\ast}-\widetilde{m}_{0})\bigr\|&\leq&\|
\mathbb{G}_{n}^{\ast}\widetilde{m}_{0}\|+\|\mathbb{G}_{n}\widetilde
{m}_{0}\|
+\|\mathbb{G}_{n}^{\ast}(\widetilde{m}_{0}-\widehat{m}^{\ast})\|
\nonumber\\
&&{}
+\|\mathbb{G}_{n}(\widetilde{m}_{0}-\widehat{m}^{\ast})\|+
\bigl\|\sqrt{n}\mathbb{P}_{n}^{\ast}\widehat{m}^{\ast}\bigr\|\\
&\equiv& L_{1}+L_{2}+L_{3}+L_{4}+L_{5}.\nonumber
\end{eqnarray}

Based on Theorem 2.2 in \cite{pw93}, we have
$L_{1}=O_{P_{W}}^{o}(1)$ in $P_{X}^{o}$-probability. The CLT
implies $L_{2}=O_{P_X}^{o}(1)$. We next consider $L_{3}$ and
$L_{4}$. By condition SB3, we can show that
$\|\widehat{\eta}^{\ast}-\eta_{0}\|=o_{P_{W}}^{o}(1)$ in
$P_{X}^{o}$-probability since $\widehat{\theta}^{\ast}$ is assumed
to be consistent, that is,
$\|\widehat{\theta}^{\ast}-\theta_{0}\|=o_{P_{W}}^{o}(1)$ in
$P_{X}^{o}$-probability, and by (\ref{term0}) and (\ref{term7}) in
Lemma \ref{term}. Then, we have $L_{3}=o_{P_{W}}^{o}(1)$ in
$P_{X}^{o}$-probability based on Lemma \ref{bmodu} and
(\ref{term7}) in Lemma \ref{term}. Next, we obtain that
$L_{4}=o_{P_{W}}^{o}(1)$ in $P_{X}^{o}$-probability based on
condition S1 and (\ref{term5}) in Lemma \ref{term}. Finally,
$L_{5}=o_{P_{\XW}}^o(1)$ based on (\ref{exactsol-b}). In summary,
(\ref{equarat}) can be rewritten as:
%
%e62 ###
%
\begin{equation}\label{equarat2}
\bigl\|\sqrt{n}P_X(\widehat{m}^{\ast}-\widetilde{m}_{0})\bigr\|\leq
O_{P_{W}}^{o}(1)+O_{P_X}^{o}(1)
\end{equation}
in $P_{X}^{o}$-probability.

Let $\alpha_{n}=\|\widehat{\theta}^{\ast}-\theta_{0}\|$. Combining
(\ref{smooth}) with (\ref{equarat2}) and noticing
(\ref{conrate-b}), we have
%
%e63 ###
%
\begin{equation}\label{rootnpf}
\sqrt{n} \|A\alpha_{n}\|\leq
O_{P_{W}}^{o}(1)+O_{P_X}^{o}(1)+O_{P_{W}}^{o}
\bigl(\sqrt{n}\alpha_{n}^{2}\vee
n^{-2\gamma+1/2} \bigr)
\end{equation}
in $P_{X}^{o}$-probability. By considering the consistency of
$\widehat{\theta}^{\ast}$ and condition I, we complete the proof
of (\ref{bconratep}) based on (\ref{rootnpf}).

We next prove (\ref{bcons}). Write
\begin{eqnarray*}
I_1&=&-\mathbb{G}_{n}^{\ast}(\widehat{m}^{\ast}-\widetilde{m}_{0})=
\sqrt{n}(\mathbb{P}_n^{\ast}-\mathbb{P}_n)(\widetilde{m}_{0}-
\widehat{m}^{\ast}),\\
I_2&=&\mathbb{G}_{n}(\widehat{m}-\widetilde{m}_{0})=
\sqrt{n}(\mathbb{P}_n-P_X)(\widehat{m}-\widetilde{m}_{0}),\\
I_3&=&-\mathbb{G}_{n}(\widehat{m}^{\ast}-\widetilde{m}_{0})=
\sqrt{n}(\mathbb{P}_n-P_X)(\widetilde{m}_{0}-\widehat{m}^{\ast}),\\
I_4&=&\sqrt{n}\mathbb{P}_{n}^{\ast}\widehat{m}^{\ast}-\sqrt
{n}\mathbb{P}_{n}
\widehat{m}.
\end{eqnarray*}
By some algebra, we obtain that
$\sqrt{n}P_X(\widehat{m}^{\ast}-\widehat{m})
+\mathbb{G}_{n}^{\ast}\widetilde{m}_{0}=\sum_{j=1}^4 I_j$.

By the definition (\ref{inter9}), we can show that $A_n\times
B_n=O_{P_W}^o(1)$ in $P_X^o$-probability if $A_n$ and $B_n$ are
both of the order $O_{P_W}^o(1)$ in $P_X^o$-probability. Then the
root-$n$ consistency of $\widehat{\theta}^{\ast}$ proven in
(\ref{bconratep}) together with SB3 implies
%
%e64 ###
%
\begin{equation}\label{inter100}
\|\widehat{\eta}^{\ast}-\eta_{0}\|\vee\|\widehat{\theta}^{\ast
}-\theta
_{0}\|
=O_{P_{W}}^{\ast}(n^{-\gamma})
\end{equation}
in $P_{X}^{o}$-probability. Thus, by Lemma \ref{bmodu}, we know
$I_{1}=O_{P_{W}}^{o}(n^{-\gamma})$ in $P_{X}^{o}$-probability.
Note that (\ref{conmod}) and (\ref{conmod2}) of condition S1
imply
%
%e65 ###
%
\begin{equation}\label{inter1}
\mathbb{G}_{n}\bigl(\widetilde{m}(\theta,\eta)-\widetilde{m}_0\bigr)
=O_{P_X}^{o}(\|\theta-\theta_{0}\|\vee\|\eta-\eta_0\|)
\end{equation}
for $(\theta,\eta)$ in the shrinking neighborhood $\mathcal{C}_n$ of
$(\theta_0,\eta_0)$.
Considering (\ref{inter1}), S3 and Proposition \ref{asythm}, we
have $I_{2}=O_{P_X}^{o}(n^{-\gamma})$. By (\ref{inter100}),
(\ref{inter1}) and (\ref{term6}), we know the order of $I_{3}$ is
$O_{P_{W}}^{o}(n^{-\gamma})$ in $P_{X}^{o}$-probability. We also obtain
$I_{4}=o_{P_X}^o(1)+o_{P_{\XW}}^o(1)$ by using (\ref{exactsol}) and
(\ref{exactsol-b}).

Therefore, we have established
%
%e66 ###
%
\begin{equation}\label{inter5}
\sqrt{n}P_X(\widehat{m}^{\ast}-\widehat{m})=-
\mathbb{G}_{n}^{\ast}\widetilde{m}_{0}+
o_{P_X}^{o}(1)+o_{P_{W}}^{o}(1)
\end{equation}
in $P_{X}^{o}$-probability. To analyze the left-hand side of
(\ref{inter5}), we rewrite it as
$\sqrt{n}P_X(\widehat{m}^{\ast}-\widetilde{m}_{0})-\sqrt{n}
P_X(\widehat{m}-\widetilde{m}_{0})$. Applying condition S2, we
obtain
%
%e67 ###
%
\begin{eqnarray}\label{for1}\quad
&&\sqrt{n}P_X\bigl(m_{11}(\theta_{0},\eta_{0})-m_{21}(\theta_{0},
\eta_{0})[H_{0}^{\dag}]\bigr)(\widehat{\theta}^{\ast}-\widehat{\theta
})\nonumber\\
&&\qquad=
-\mathbb{G}_{n}^{\ast}\widetilde{m}_{0}+o_{P_X}^{o}(1)
+o_{P_{W}}^{o}(1)+O_{P_X}^{o}(n^{1/2-2\gamma})+
O_{P_{W}}^{o}(n^{1/2-2\gamma})\\
&&\qquad=
-\mathbb{G}_{n}^{\ast}\widetilde{m}_{0} +o_{P_X}^{o}(1)+
o_{P_{W}}^{o}(1)\nonumber
\end{eqnarray}
in $P_{X}^{o}$-probability, by considering condition S3, SB3 and
the range of $\gamma$. Note that $o_{P_X}^o(1)$ in (\ref{for1}) is also
of the order $o_{P_{\XW}}^o(1)$, and thus is of the order $o_{P_W}^o(1)$
in $P_X^o$-probability by (\ref{term0}). Moreover,
according to condition I we have that
$A=P_X(m_{11}(\theta_{0},\eta_{0})-m_{21}(\theta_{0},
\eta_{0})[H_{0}^{\dag}])$ is nonsingular. We obtain (\ref{bcons})
by multiplying $A^{-1}$ on both sides of (\ref{for1}).

By applying Lemma 4.6 in \cite{pw93} under the bootstrap weight
conditions, we obtain (\ref{proconv}).
Proposition \ref{asythm} together with Lemma 2.11 in \cite{v98}
implies that
%
%e68 ###
%
\begin{equation}\label{proconv2}
\sup_{x\in\mathbb{R}^d} \bigl|P_X\bigl(\sqrt{n}(\widehat{\theta
}-\theta_0)\leq
x\bigr) -P\bigl(N(0,\Sigma)\leq x\bigr) \bigr|=o(1).
\end{equation}
Combining (\ref{proconv}) and (\ref{proconv2}), we obtain
(\ref{probc1}).

\begin{appendix}\label{app}
\section*{Appendix}

%s7.1 ###
\subsection{Measurability and stochastic orders}\label{appA1}

\textit{Measurability condition} $M(P)$: we say that a class of
random functions $\mathcal{F}\in M(P)$ if $\mathcal{F}$ is nearly
linearly deviation measurable for $P$ and that both
$\mathcal{F}^{2}$ and $\mathcal{F}'^{2}$ are nearly linearly
supremum measurable for $P$. Here $\mathcal{F}^{2}$ and
$\mathcal{F}'^{2}$ denote the classes of squared functions and
squared differences of functions from $\mathcal{F}$, respectively.
It is known that if $\mathcal{F}$ is countable, or if
$\{\mathbb{P}_{n}\}_{n=1}^{\infty}$ are stochastically separable
in $\mathcal{F}$, or if $\mathcal{F}$ is image admissible
Suslin~\cite{d84}, then $\mathcal{F}\in M(P)$. More precise descriptions
can be found in pages 853 and 854 of \cite{gz90}.

The following lemma is very important since it accurately
describes the transition of stochastic orders among different
probability spaces. We implicitly assume the random quantities in
Lemma \ref{term} posses enough measurability so that the usual
Fubini theorem can be used freely.
\begin{lemma}\label{term}
Suppose that
\begin{eqnarray*}
Q_n&=&o_{P_{W}}^{o}(1) \qquad\mbox{in }
P_{X}^{o}\mbox{-probability},\\
R_{n}&=&O_{P_{W}}^{o}(1) \qquad\mbox{in }
P_{X}^{o}\mbox{-probability}.
\end{eqnarray*}
We have
%
%e73 ###
%e72 ###
%e71 ###
%e70 ###
%e69 ###
%
\begin{eqnarray}\label{term0}
A_{n}=o_{P_{\XW}}^{o}(1) \quad\Longleftrightarrow\quad
A_{n}&=&o_{P_{W}}^{o}(1) \qquad\mbox{in } P_{X}^{o}
\mbox{-probability},\\
\label{term1}
B_{n}=O_{P_{\XW}}^{o}(1) \quad\Longleftrightarrow\quad
B_{n}&=&O_{P_{W}}^{o}(1) \qquad\mbox{in }
P_{X}^{o}\mbox{-probability},\\
\label{term5}
C_n=Q_n\times O_{P_X}^{o}(1) \quad\Longrightarrow\quad\hspace*{1pt}
C_{n}&=&o_{P_{W}}^{o}(1) \qquad\mbox{in }
P_{X}^{o}\mbox{-probability},\\
\label{term6}
D_n=R_n\times O_{P_X}^{o}(1) \quad\Longrightarrow\quad
D_{n}&=&O_{P_{W}}^{o}(1) \qquad\mbox{in }
P_{X}^{o}\mbox{-probability},\\
\label{term7}
E_n=Q_n\times
R_n \quad\Longrightarrow\quad E_n&=&o_{P_W}^o(1) \qquad\mbox{in }
P_{X}^{o}\mbox{-probability}.
\end{eqnarray}
\end{lemma}
\begin{pf}
To verify (\ref{term0}), we have for every $\varepsilon,
\nu>0$,
%
%e74 ###
%
\begin{eqnarray}\label{inetec1}
P_{X}^{o} \{P_{W|X}^{o}(|A_{n}|\geq\varepsilon)\geq\nu
\}&\leq&\frac{1}{\nu}E_X^{o}P_{W|X}^{o}(|A_{n}|\geq\varepsilon)
\nonumber\\[-8pt]\\[-8pt]
&\leq&
\frac{1}{\nu}E^{o}_XE_{W|X}^{o}1\{|A_{n}|\geq\varepsilon\}\nonumber
\end{eqnarray}
by Markov's inequality. According to Lemmas 6.5 and 6.14 in
\cite{k08a}, we have $E_X^{o}E_{W|X}^{o}1\{|A_{n}|\geq\varepsilon\}
\leq E_{\XW}^{o}1\{|A_{n}|\geq\varepsilon\}=P_{\XW}^{o}
(|A_{n}|\geq\varepsilon)$, and thus
%
%e75 ###
%
\begin{equation}\label{inter2}
P_{X}^{o} \{P_{W|X}^{o}(|A_{n}|\geq\varepsilon)\geq\nu
\}\leq\frac{1}{\nu}P_{\XW}^o(|A_n|\geq\varepsilon).
\end{equation}
From (\ref{inter2}), we\vspace*{1pt} can conclude that if
$A_{n}=o_{P_{\XW}}^{o}(1)$, then $A_{n}=o_{P_W}^o(1)$ in
$P_{X}^{o}$-probability. Another direction of (\ref{term0})
follows from the following inequalities: for any $\varepsilon,
\eta>0$,
%
%e76 ###
%
\begin{eqnarray}\label{inetec2}
P_{\XW}^{o}(|A_n|\geq\varepsilon)&=&E^{o}_X\{P_{W|X}^{o}(|A_{n}|
\geq\varepsilon)\}\nonumber\\
&=&E^{o}_X\bigl\{P_{W|X}^{o}(|A_n|\geq\varepsilon)1\{P_{W|X}^{o
}(|A_n|\geq\varepsilon)\geq\eta\}\bigr\}\nonumber\\
&&{} +
E^{o}_X\bigl\{P_{W|X}^{o}(|A_n|\geq\varepsilon)1\{P_{W|X}^{o
}(|A_n|\geq\varepsilon)<\eta\}\bigr\}\\
&\leq&E^{o}_X\bigl\{1\{P_{W|X}^{o
}(|A_n|\geq\varepsilon)\geq\eta\}\bigr\}+\eta\nonumber\\
&\leq&P_{X}^{o}\{P_{W}^{o}(|A_n|\geq
\varepsilon)\geq\eta\}+\eta.\nonumber
\end{eqnarray}
Note that the first term in (\ref{inetec2}) can be made
arbitrarily small by the assumption that $A_n=o_{P_W}^{o}(1)$ in
$P_{X}^{o}$-probability. Since $\eta$ can be chosen arbitrarily
small, we can show
$\lim_{n\rightarrow\infty}P_{\XW}^{o}(|A_n|\geq\varepsilon)=0$ for any
$\varepsilon>0$. This completes the proof of (\ref{term0}).
(\ref{term1}) can be shown similarly by using the inequalities
(\ref{inetec1}) and (\ref{inetec2}).

As for (\ref{term5}), we establish the following inequalities:
\begin{eqnarray*}
&&P_{X}^{o} \bigl\{P_{W|X}^{o}\bigl(|Q_n\times
O_{P_X}^{o}(1)|\geq\varepsilon\bigr)\geq\nu
\bigr\}\\
&&\qquad\leq P_{X}^{o} \bigl\{P_{W|X}^{o}\bigl(|Q_n|\geq\varepsilon/
|O_{P_X}^{o}(1)|\bigr)\geq\nu\bigr\}\\
&&\qquad\leq P_{X}^{o}\bigl\{P_{W|X}^{o}(|Q_n|\geq\varepsilon/M)+
P_{W|X}^{o}\bigl(|O_{P_X}^{o}(1)|\geq M\bigr)\geq\nu\bigr\}\\
&&\qquad\leq P_{X}^{o}\{P_{W|X}^{o}(|Q_n|\geq\varepsilon/M)\geq
\nu/2\}+P_{X}^{o}\bigl\{P_{W|X}^{o}\bigl(|O_{P_X}^{o}(1)|\geq
M\bigr)\geq\nu/2\bigr\}\\
&&\qquad\leq P_{X}^{o}\{P_{W|X}^{o}(|Q_n|\geq\varepsilon/M)\geq
\nu/2\}+\frac{2}{\nu}P_{X}^{o}\bigl(|O_{P_X}^{o}(1)|\geq M\bigr)
\end{eqnarray*}
for any $\varepsilon,\nu,M>0$. Since $M$ can be chosen arbitrarily
large, we can show (\ref{term5}) by considering the definition of
$O_{P_X}^{o}(1)$. The proof of (\ref{term6}) is similar by using
the above set of inequalities.
%we will establish the following inequalities:
%&&P_{X}^{o}\{P_{W|X}^{o}(|R_n\times O_{P_{X}}^{o}(1)|\geq M)\leq
%M/|O_{P_{X}}^{o}(1)|)\leq\varepsilon\}\\
%&\geq& P_{X}^{o}\{P_{W|X}^{o}(|R_{n}|\geq M/T)+
%P_{W|X}^{o}(|O_{P_{X}}^{o}(1)|\geq T)\leq\varepsilon\}\\
%&\geq&P_{X}^{o}\{P_{W|X}^{o}(|R_{n}|\geq M/T)\leq\varepsilon/2 \&
%P_{W|X}^{o}(|O_{P_{X}}^{o}(1)|\geq T)\leq\varepsilon/2\}
%for arbitrary $M,T,\varepsilon>0$. By the assumption on $R_n$, we
%know $P_{X}^{o}\{P_{W|X}^{o}(|R_n|\geq
%M/T)\leq\varepsilon/2\}\rightarrow1$ by choosing proper $M$ and $T$.
%Next, we establish
%P_{X}^{o}\{P_{W|X}^{o}(|O_{P_X}^{o}(1)|\geq T)\leq
%{\varepsilon/2}\\&\geq& 1-\frac{E_{W}^{o}P_{X}^{o}
%(|O_{P_X}^{o}(1)|\geq T)} {\varepsilon/2}.
%By choosing large enough $T$, we have shown (\ref{term6}).
The proof of (\ref{term5}) can be carried over to prove
(\ref{term7}). Similarly, we establish the following inequalities:
\begin{eqnarray*}
&&P_{X}^{o} \{P_{W|X}^{o}(|Q_n\times R_n|\geq\varepsilon)\geq\eta
\}\\
&&\qquad\leq P_{X}^{o}\{P_{W|X}^{o}(|Q_n|\geq\varepsilon/M)\geq
\eta/2\}+P_{X}^{o}\{P_{W|X}^{o}(|R_n|\geq M)\geq\eta/2\}
\end{eqnarray*}
for any $\varepsilon,\eta,M>0$. Then by selecting sufficiently large
$M$, we can show that
\[
P_{X}^{o} \{P_{W|X}^{o}(|Q_n\times
R_n|\geq\varepsilon)\geq\eta\}\rightarrow0
\]
as
$n\rightarrow\infty$ for any $\varepsilon,\eta>0$.
\end{pf}

%s7.2 ###
\subsection{Two useful inequalities}\label{appA2}

Here we give two key inequalities used in proving
Lemmas \ref{lerate} and \ref{bmodu}.

\subsubsection*{\texorpdfstring{Multiplier inequality
(Lemma 4.1 of \protect\cite{wz96})}{Multiplier inequality (Lemma 4.1 of [41])}}

\mbox{}

Let $W_{n}=(W_{n1},\ldots,W_{nn})'$ be nonnegative exchangeable
random variables on $(\mathcal{W},\Omega,P_{W})$ such that, for
every $n$, $R_{n}=\int_{0}^{\infty}\sqrt{P_{W}(W_{n1}\geq
u)}\,du<\infty$. Let $Z_{ni}$, $i=1,2,\ldots,n$, be i.i.d. random
elements in
$(\mathcal{X}^{\infty},\mathcal{A}^{\infty},P_X^{\infty})$ with
values in $\ell^{\infty}(\mathcal{F}_{n})$, and write
$\|\cdot\|_{n}={\sup_{f\in\mathcal{F}_{n}}}|Z_{ni}(f)|$. It is
assumed that $Z_{ni}$'s are independent of $W_{n}$. Then for any
$n_{0}$ such that $1\leq n_{0}<\infty$ and any $n>n_{0}$, the
following inequality holds:
%
%e77 ###
%
\begin{eqnarray}\label{mulineq}
E_{\XW}^{o} \Biggl\|\frac{1}{\sqrt{n}}\sum_{i=1}^{n}
W_{ni}Z_{ni} \Biggr\|_{n}&\leq&
n_{0}E^{o}_X\|Z_{n1}\|_{n}\cdot\frac{E_{W}(\max_{1\leq i\leq
n}W_{ni})}{\sqrt{n}}\nonumber\\[-8pt]\\[-8pt]
&&{} + R_{n}\cdot\max_{n_{0}< i\leq
n} \Biggl\{E^{o}_X\frac{1}{\sqrt{i}}
\Biggl\|\sum_{j=n_{0}+1}^{i}Z_{nj} \Biggr\|_{n} \Biggr\}.\nonumber
\end{eqnarray}

\subsubsection*{\texorpdfstring{Hoffmann--Jorgensen inequality for moments (Proposition
\textup{A.1.5} in \protect\cite{vw96})}{Hoffmann--Jorgensen inequality for moments (Proposition
\textup{A.1.5} in [38])}}

\mbox{}

Let $1\leq p<\infty$ and suppose that $V_{1},\ldots,V_{n}$ are
independent stochastic processes with mean zero indexed by an
arbitrary index set $T$. Then there exist constants $K_{p}$ and
$0<v_{p}<1$ such that
\[
E^{o} \Biggl\|\sum_{i=1}^{n}V_{i} \Biggr\|^{p}\leq
K_{p} \Bigl\{{E^{o}\max_{1\leq k\leq
n}}\|V_{k}\|^{p}+[G^{-1}(v_{p})]^{p} \Bigr\},
\]
where $\|Y\|=\sup_{t}|Y_{t}|$ denotes the supremum of a stochastic process
$\{Y_{t}, t\in T\}$, and
$G^{-1}(v)=\inf\{u\dvtx P^o (\|\sum_{i=1}^{n}V_{i}\|\leq
v )\geq u \}$.

%s7.3 ###
\subsection{\texorpdfstring{Proof of Lemma \protect\ref{bmodu}}{Proof of Lemma 2}}\label{appA3}

We first write
$\mathbb{G}_{n}^{\ast}(\widetilde{m}(\theta,\eta)-\widetilde{m}_{0})$
as the sum of
$\mathbb{G}_{n}^{\ast}(\widetilde{m}(\theta,\eta)-\widetilde
{m}(\theta_
{0},\eta))$ and
$\mathbb{G}_{n}^{\ast}(\widetilde{m}(\theta_{0},\eta)-\widetilde{m}_{0})$.
By the Taylor expansion, the first term becomes
$(\theta-\theta_0)'\mathbb{G}_n^{\ast}(\partial/\partial\theta)
\widetilde{m}(\bar{\theta},\eta)$, where $\bar{\theta}$ is
between $\theta$ and $\theta_0$. By SB2 and Theorem 2.2 in
\cite{pw93}, we know that the first term is of the order
$O_{P_{W}}^{o}(\|\theta-\theta_{0}\|)$ in $P_{X}^{o}$-probability.
We next consider the second term. Let
%
%e78 ###
%
\begin{equation}\label{deltadfn}
\Delta_{n}=\sup_{\eta\in
U_{n}} \biggl\{\frac{\|\mathbb{G}_{n}^{\ast}(\widetilde{m}(\theta
_{0},\eta
)-\widetilde{m}_{0})\|}
{\|\eta-\eta_{0}\|} \biggr\},
\end{equation}
where $U_{n}=\{\eta\dvtx\|\eta-\eta_{0}\|\leq\delta_{n}\}$ for any
$\delta_{n}\rightarrow0$. Note that we can write
$\Delta_{n}=\|\mathbb{G}_{n}^{\ast}\|_{\mathcal{S}_{n}}$, where
$\|\mathbb{G}_{n}^{\ast}\|_{\mathcal{S}_{n}}={\sup_{f\in\mathcal{S}_{n}
}}|\mathbb{G}_{n}^{\ast}f|$. By (\ref{term1}), to verify the
bootstrap equicontinuity condition that
$\mathbb{G}_{n}^{\ast}(\widetilde{m}(\theta_{0},\eta)-
\widetilde{m}_{0})=O_{P_{W}}^{o}(\|\eta-\eta_{0}\|)$ in
$P_X^o$-probability, it suffices to show
%
%e79 ###
%
\begin{equation}\label{bounddel}
\lim\sup_{n\rightarrow\infty}E_{\XW}^{o}\Delta_{n}< \infty.
\end{equation}

Note that
\[
\mathbb{G}_{n}^{\ast}=\frac{1}{\sqrt{n}}\sum
_{i=1}^{n}(W_{ni}-1)\delta
_{X_{i}}=\frac{1}{\sqrt{n}}\sum_{i=1}^{n}
(W_{ni}-1)(\delta_{X_{i}}-P_X)
\]
by condition W2. Let $W_{n}'=(W_{n1}',\ldots,W_{nn}')$ be
exchangeable bootstrap\break weights generated from $P_{W'}$, an
independent copy of $P_{W}$. The bootstrap weight conditions W1
and W2 imply that $E_{W'}W_{ni}'=1$ for $i=1,\ldots,n$. Let
\[
m_{n}(\eta,\eta_0)=\frac{\widetilde{m}(\theta_{0},\eta)-
\widetilde{m}_{0}}{\|\eta-\eta_{0}\|}.
\]
Then we have
\begin{eqnarray*}
E_{\XW}^{o}\Delta_{n}&=&{E_{\XW}^{o}\sup_{\eta\in
U_{n}}}\|\mathbb{G}_{n}^{\ast}m_{n}(\eta,\eta_{0})\|\\&=&
E_{\XW}^{o}\sup_{\eta\in
U_{n}} \Biggl\|\frac{1}{\sqrt{n}}\sum_{i=1}^{n}(W_{ni}-1)(\delta_{X_{i}}
-P_X)m_{n}(\eta,\eta_{0}) \Biggr\|\\
&=&E_{\XW}^{o}\sup_{\eta\in
U_{n}} \Biggl\|\frac{1}{\sqrt{n}}\sum_{i=1}^{n}(W_{ni}-E_{W'}
W_{ni}')(\delta_{X_{i}}-P_X)m_{n}(\eta,\eta_{0})
\Biggr\|\\
&\leq&E_{\XW}^{o}E^{o}_{W'}\sup_{\eta\in
U_{n}} \Biggl\|\frac{1}{\sqrt{n}}\sum_{i=1}^{n}(W_{ni}-W_{ni}')
(\delta_{X_{i}}-P_X)m_{n}(\eta,\eta_{0}) \Biggr\|.
\end{eqnarray*}
%
%The last inequality follows from the Jensen's inequality since
%$\|\cdot\|_{\mathcal{S}_{n}}$ is convex.
To further bound $E_{\XW}^{o}\Delta_{n}$, we employ the
symmetrization argument familiar in the empirical process
literature to obtain
%
%e80 ###
%
\begin{eqnarray}\label{inter200}
E_{\XW}^{o}\Delta_{n}&\leq&E_{\XW}^{o} \sup_{\eta\in
U_{n}} \Biggl\|\frac{1}{\sqrt{n}}\sum_{i=1}^{n}W_{ni}
(\delta_{X_{i}}-P_X)m_{n}(\eta,\eta_{0}) \Biggr\|\nonumber\\
&&{} +
E_{\XW}^{o}E^{o}_{W'}\sup_{\eta\in
U_{n}} \Biggl\|\frac{1}{\sqrt{n}}\sum_{i=1}^{n}W_{ni}'
(\delta_{X_{i}}-P_X)m_{n}(\eta,\eta_{0}) \Biggr\|\\
&=& 2E_{\XW}^{o}\sup_{\eta\in
U_{n}} \Biggl\|\frac{1}{\sqrt{n}}\sum_{i=1}^{n}W_{ni}
(\delta_{X_{i}}-P_X)m_{n}(\eta,\eta_{0}) \Biggr\|.\nonumber
\end{eqnarray}

We next apply the multiplier inequality (\ref{mulineq}) to (\ref{inter200})
with $Z_{ni}=\{(\delta_{X_{i}}-P_X)m_{n}(\eta,\eta_{0})\dvtx\eta\in
U_{n}\}$.
Define
\[
\|Z_{ni}\|_{n}={\sup_{\eta\in
U_{n}}} \|(\delta_{X_{i}}-P_X)m_{n}(\eta,\eta_{0}) \|.
\]
To show (\ref{bounddel}), we need only to show
%
%e81 ###
%
\begin{equation}\label{inter20}
E_{W}\Bigl(\max_{1\leq i\leq
n}W_{ni}\Bigr)\big/\sqrt{n}\rightarrow0,
\end{equation}
$\lim\sup_n E_X^o\|Z_{n1}\|_n<\infty$, and
%
%e82 ###
%
\begin{equation}\label{lastterm}
\lim\sup_n \max_{n_{0}< i\leq n}E^{o}_X \sup_{\eta\in
U_{n}} \Biggl\|\frac{1}{\sqrt{i}}\sum_{j=n_{0}+1}^{i}
Z_{ni} \Biggr\|<\infty
\end{equation}
for some $n_0< \infty$.
The bootstrap weight conditions W3 and W4 together with Lemma 4.7 in
\cite{pw93} imply (\ref{inter20}).
Note that
\begin{eqnarray*}
E^{o}_X\|Z_{n1}\|_{n}&=&{E^{o}_{X}\sup_{\eta\in
U_{n}}}\|(\delta_{X_{1}}-P_X)m_{n}(\eta,\eta_{0})\|\\
&\leq&{E^{o}_X\sup_{\eta\in
U_{n}}}\|m_{n}(\eta,\eta_{0})(X_{1})\|+{E_X^o\sup_{\eta\in U_n}}\|E_X
m_{n}(\eta,\eta_0)\|\\
&\leq&2E_X^{o}S_{n}(X_{1}),
\end{eqnarray*}
where $S_{n}$ is the envelop of the class $\mathcal{S}_{n}$
defined in (\ref{sndfn}), and the first inequality follows from
the Fatou's lemma. Condition SB1 implies
%
%e84 ###
%e83 ###
%
\begin{eqnarray}
\label{addimp1}
\frac{1}{\sqrt{n}}E^{o}_X\max_{1\leq k\leq
n}S_{n}(X_{k})&\longrightarrow&
0,\\
\label{addimp2}
\lim\sup_{n\rightarrow\infty}E^{o}_XS_{n}(X_{1})&<&\infty;
\end{eqnarray}
see page 120 of \cite{vw96}. The result (\ref{addimp2}) implies $\lim
\sup_{n}
E^{o}_X\|Z_{n1}\|_{n}<\infty$.

It remains to show (\ref{lastterm}).
We apply the Hoffmann--Jorgensen
inequality with \mbox{$p=1$} in Appendix \ref{appA2}. First, we
establish
%
%e85 ###
%
\begin{eqnarray}\label{lastterm1}
E^{o}_X\sup_{\eta\in
U_{n}} \Biggl\|\frac{1}{\sqrt{n}}\sum_{i=1}^{n}Z_{ni} \Biggr\|
&\leq&
K_{1} \biggl\{\frac{1}{\sqrt{n}}{E^{o}_X\max_{1\leq k\leq
n}}\|Z_{nk}\|_{n}+G_{n}^{-1}(v_{1}) \biggr\}\nonumber\\[-8pt]\\[-8pt]
&\leq&
I_{1}+I_{2},\nonumber
\end{eqnarray}
where $K_{1}$ and $0<v_{1}<1$ are constants and
\[
G_{n}(t)=P_{X}^{o} \Biggl(n^{-1/2} \Biggl\|\sum_{i=1}^{n}Z_{ni} \Biggr\|
_{n}\leq
t \Biggr).
\]
Obviously, (\ref{addimp1}) implies that $I_{1}\rightarrow0$. We
next consider $I_{2}$. Note that assumption S1 implies
$\|\mathbb{G}_{n}\|_{\mathcal{S}_{n}}=\|n^{-1/2}\sum_{i=1}^{n}
Z_{ni}\|_{n}=O_{P_{X}}^{o}(1)$. Hence, there
exists a finite constant $M_{t}$ such that
$\lim\inf_{n}G_{n}(M_{t})\geq t$ for every $1>t>0$.
It follows that
$\lim\sup_{n}G_{n}^{-1}(v_{1})\leq M_{v_{1}}<\infty$ since
$0<v_{1}<1$. Thus, the left-hand side of (\ref{lastterm1})
is bounded away from infinity, and therefore
(\ref{lastterm}) holds in light of the following result
from the triangular inequality
\begin{eqnarray*}
\max_{n_{0}< i\leq n}E^{o}_X \sup_{\eta\in
U_{n}} \Biggl\|\frac{1}{\sqrt{i}}\sum_{j=n_{0}+1}^{i}
Z_{nj} \Biggr\|
&\leq&\max_{n_{0}< i\leq n}E^{o}_X \sup_{\eta\in
U_{n}} \Biggl\|\frac{1}{\sqrt{i}}\sum_{j=1}^{i}
Z_{nj} \Biggr\|\\
&&{}+E_X^o\sup_{\eta\in U_n} \Biggl\|\frac{1}{\sqrt
{n_0}}\sum
_{j=1}^{n_0}Z_{nj} \Biggr\|.
\end{eqnarray*}
The proof of Lemma \ref{bmodu} is complete.

%s7.4 ###
\subsection{\texorpdfstring{Proof of Theorem \protect\ref{ratethm2}}{Proof of Theorem 2}}\label{appA4}

Using (\ref{neunb-b}) and the fact that $U(\theta_0,\eta_0)=0$, we have
%
%e86 ###
%
\begin{eqnarray}
&&
U(\widetilde{\theta},\widehat{\eta}^{\ast}_{\widetilde{\theta}})
-U(\theta_{0},\eta_{0})\nonumber\\
&&\qquad=
U(\widetilde{\theta},\widehat{\eta}^{\ast}_{\widetilde{\theta}})
-U_{n}^{\ast}(\widetilde{\theta},\widehat{\eta}_{\widetilde
{\theta
}}^{\ast})+
O_{P_{\XW}}^{o}(n^{-1/2})\nonumber\\[-8pt]\\[-8pt]
&&\qquad=-(U_{n}^{\ast}-U_{n})(\widetilde{\theta},\widehat{\eta}
_{\widetilde{\theta}}^{\ast})-
(U_{n}-U)(\widetilde{\theta},\widehat{\eta}_{\widetilde{\theta}}
^{\ast})+
O_{P_{\XW}}^{o}(n^{-1/2})\nonumber\\
&&\qquad=L_1+L_2+O_{P_{\XW}}^{o}(n^{-1/2}).\nonumber
\end{eqnarray}
Further, based\vspace*{1pt} on conditions (\ref{ratecon1}) and
(\ref{ratecon7}), we apply Lemma 4.2 in \cite{wz96} to obtain
that
$L_1=-(U_{n}^{\ast}-U_{n})(\theta_0,\eta
_0)+o_{P_{\XW}}^o(n^{-1/2}\vee
\|\widetilde{\theta}-\theta_0\|\vee\|\widehat{\eta}_{\widetilde
{\theta}}
^{\ast}-\eta_0\|)$. By Lemma 3.3.5 in \cite{vw96} given
(\ref{ratecon1}) and (\ref{ratecon2}), we have
$L_2=-(U_{n}-U)(\theta_0,\eta_0)+o_{P_{\XW}}^o(n^{-1/2}\vee
\|\widetilde{\theta}-\theta_0\|\vee\|\widehat{\eta}_{\widetilde
{\theta}}
^{\ast}-\eta_0\|)$. By applying CLT and Theorem 2.2 in \cite{pw93}
under condition (\ref{ratecon1}) to $L_1$ and $L_2$, we have
%
%e87 ###
%
\begin{equation}\label{rhs1}\quad
U(\widetilde{\theta},\widehat{\eta}^{\ast}_{\widetilde{\theta}})
-U(\theta_{0},\eta_{0})=O_{P_{\XW}}^{o}(n^{-1/2})+o_{P_{\XW}}^{o}
(\|\widetilde{\theta}-\theta_{0}\|\vee\|\widehat{\eta
}_{\widetilde
{\theta}}^{\ast}-\eta_{0}\|).
\end{equation}

We next apply the Taylor expansion to
get
\begin{eqnarray*}
&&U(\widetilde{\theta},\widehat{\eta}^{\ast}_{\widetilde{\theta}})-
U(\theta_{0},\eta_{0})\\
&&\qquad=\dot{U}(\widetilde{\theta}-\theta_{0},\widehat{\eta
}_{\widetilde
{\theta}}
^{\ast}-\eta_{0})+o(\|\widetilde{\theta}-\theta_{0}\|
\vee\|\widehat{\eta}_{\widetilde{\theta}}^{\ast}-\eta_{0}\|)\\
&&\qquad=
\dot{U}(\widetilde{\theta}-\theta_{0},0)+
\dot{U}(0,\widehat{\eta}_{\widetilde{\theta}}^{\ast}-\eta_{0})
+o(\|\widetilde{\theta}-\theta_{0}\|
\vee\|\widehat{\eta}_{\widetilde{\theta}}^{\ast}-\eta_{0}\|)
\end{eqnarray*}
by the assumed Fr\'{e}chet differentiability of $U$ and linearity
of $\dot{U}$. Note that $U$ has bounded Fr\'{e}chet derivative and
$\dot{U}(0,\cdot)$ is continuously invertible. Thus, we can
conclude that
\[
U(\widetilde{\theta},\widehat{\eta}^{\ast}_{\widetilde{\theta}})-
U(\theta_{0},\eta_{0})\geq c\|\widehat{\eta}_
{\widetilde{\theta}}^{\ast}
-\eta_{0}\|+O(\|\widetilde{\theta}-\theta_{0}\|)+o(\|\widetilde
{\theta
}-\theta_{0}\|
\vee\|\widehat{\eta}_{\widetilde{\theta}}^{\ast}-\eta_{0}\|)
\]
for some $c>0$. Combining the above inequality with (\ref{rhs1}),
we can establish the following inequality:
\[
\|\widehat{\eta}_{\widetilde{\theta}}^{\ast}-\eta_{0}\|\lesssim
O_{P_{\XW}}^{o}(\|\widetilde{\theta}-\theta_{0}\|\vee
n^{-1/2})+o_{P_{\XW}}^{o}
(\|\widehat{\eta}_{\widetilde{\theta}}^{\ast}-\eta_{0}\|),
\]
which implies (\ref{nuirate}).

%s7.5 ###
\subsection{\texorpdfstring{Proof of Theorem \protect\ref{ratethm}}{Proof of Theorem 3}}
\label{appA5}

According to (\ref{term1}), we need only to show that
%
%e88 ###
%
\begin{equation}\label{ratedes}
P_{\XW}^{o} \bigl(d(\widehat{\eta}^{\ast}_{\widetilde{\theta}},\eta
_{n})\geq
2^{M_{n}}(\delta_{n}\vee\|\widetilde{\theta}-\theta_{0}\|),
\widetilde{\theta}\in\Theta,
\widehat{\eta}_{\widetilde{\theta}}^{\ast}\in\mathcal{H}_{n}
\bigr)\longrightarrow0
\end{equation}
as $n\rightarrow\infty$ and $M_{n}\rightarrow\infty$. The basic
idea in
proving (\ref{ratedes}) is first
to partition the whole parameter space into ``shells,'' and then bound
the probability of each shell under conditions
(\ref{ratecon3})--(\ref{ratecon5}).

For now we fix $M=M_{n}$ and then allow it to increase to
infinity. We first define the shell $S_{n,j,M}$ as
\[
S_{n,j,M}=\{(\theta,\eta)\in\Theta\times\mathcal
{H}_{n}\dvtx2^{j-1}\delta_{n}<
d(\eta,\eta_{n})\leq2^{j}\delta_{n},d(\eta,\eta_{n})\geq2^{M}
\|\theta-\theta_{0}\|\}
\]
with $j$ ranging over the integers and $M>0$. Obviously, the event
$\{\widetilde{\theta}\in\Theta,
\widehat{\eta}_{\widetilde{\theta}}^{\ast}\in\mathcal
{H}_{n}\dvtx d(\widehat
{\eta}^{\ast}_{\widetilde{\theta}},\eta_{n})\geq
2^{M}(\delta_{n}\vee\|\widetilde{\theta}-\theta_{0}\|)\}$ is
contained in the union of the events
$\{(\widetilde{\theta},\widehat{\eta}_{\widetilde{\theta}}^{\ast
})\in
S_{n,j,M}\}$ for $j\geq M$. Thus, we have
\begin{eqnarray*}
&&
P_{\XW}^{o} \bigl(d(\widehat{\eta}^{\ast}_{\widetilde{\theta}},\eta
_{n})\geq
2^{M}(\delta_{n}\vee\|\widetilde{\theta}-\theta_{0}\|),
\widetilde{\theta}\in\Theta,
\widehat{\eta}_{\widetilde{\theta}}^{\ast}\in\mathcal
{H}_{n} \bigr)\\
&&\qquad\leq
\sum_{j\geq
M}P_{\XW}^{o} \bigl((\widetilde{\theta},\widehat{\eta}_{\widetilde
{\theta
}}^{\ast})\in
S_{n,j,M} \bigr)\\
&&\qquad\leq\sum_{j\geq M}P_{\XW}^{o} \Bigl(\sup_{(\theta,\eta)\in
S_{n,j,M}}\mathbb{P}_{n}^{\ast}\bigl(v(\theta,\eta)-v(\theta,\eta
_{n})\bigr)\geq
0 \Bigr).
\end{eqnarray*}
The second inequality follows from the definition of
$\widehat{\eta}^{\ast}_{\widetilde{\theta}}$. By the smoothness
condition on $v(\theta,\eta)$, that is, (\ref{ratecon3}), we have the
following inequality when $(\theta,\eta)\in S_{j,n,M}$ for $j\geq M$:
%
%e89 ###
%
\begin{equation}\label{smooth2}
P_X\bigl(v(\theta,\eta)-v(\theta,\eta_{n})\bigr)\lesssim-d(\eta,\eta
_{n})^{2}+\|
\theta-\theta_{0}\|^{2}\lesssim
-2^{2j-2}\delta_{n}^{2}
\end{equation}
for sufficiently large $M$.

Considering (\ref{smooth2}), we have
\begin{eqnarray*}
&&P_{\XW}^{o} \bigl(d(\widehat{\eta}^{\ast}_{\widetilde{\theta
}},\eta
_{n})\geq
2^{M}(\delta_{n}\vee\|\widetilde{\theta}-\theta_{0}\|),
\widetilde{\theta}\in\Theta,
\widehat{\eta}_{\widetilde{\theta}}^{\ast}\in\mathcal
{H}_{n} \bigr)\\
&&\qquad\leq\sum_{j\geq M}P_{\XW}^{o} \Bigl(\sup_{(\theta,\eta)\in
S_{n,j,M}}\sqrt{n}(\mathbb{P}_{n}^{\ast}-P_X)\bigl(v(\theta,\eta
)-v(\theta
,\eta_{n})\bigr)\gtrsim\sqrt{n}2^{2j-2}\delta_{n}^{2} \Bigr)\\
&&\qquad\leq\sum_{j\geq M}P_{\XW}^{o} \Bigl(\sup_{(\theta,\eta)\in
S_{n,j,M}}\bigl|\mathbb{G}_{n}^{\ast}\bigl(v(\theta,\eta)-v(\theta,\eta
_{n})\bigr)\bigr|\gtrsim\sqrt{n}2^{2j-3}
\delta_{n}^{2} \Bigr)\\
&&\qquad\quad{} +P_{X}^{o} \Bigl( \sup_{(\theta,\eta)\in
S_{n,j,M}}\bigl|\mathbb{G}_{n}\bigl(v(\theta,\eta)-v(\theta,\eta
_{n})\bigr)\bigr|\gtrsim\sqrt
{n}2^{2j-3}\delta_{n}^{2} \Bigr)\\
&&\qquad\lesssim\sum_{j\geq
M}\frac{\psi_{n}^{\ast}(2^{j}\delta_{n})}{\sqrt{n}\delta_{n}^{2}2^{2j}}+
\frac{\psi_{n}(2^{j}\delta_{n})}{\sqrt{n}\delta_{n}^{2}2^{2j}}\\
&&\qquad\lesssim\sum_{j\geq M}2^{j(\alpha-2)},
\end{eqnarray*}
where the third inequality follows from the Markov inequality and
(\ref{ratecon4}) and (\ref{ratecon5}). Note that the assumption that
$\delta\mapsto\psi_{n}(\delta)/\delta^{\alpha}$
[$\delta\mapsto\psi_{n}^{\ast}(\delta)/\delta^{\alpha}$] is
decreasing for some $0<\alpha<2$ implies that
$\psi_{n}(c\delta)\leq c^{\alpha}\psi_{n}(\delta)$ for every
$c>1$. Combining these with the assumption that
$\psi_{n}(\delta_{n})\leq\sqrt{n}\delta_{n}^{2}$ and
$\psi_{n}^{\ast}(\delta_{n})\leq\sqrt{n}\delta_{n}^{2}$, we
obtain the
last inequality in the above display. By letting $M=M_{n}\rightarrow
\infty$, we complete the proof of
(\ref{ratedes}), and thus Theorem \ref{ratethm}.

%s7.6 ###
\subsection{\texorpdfstring{Proof of Lemma \protect\ref{lerate}}{Proof of Lemma 1}}\label{appA6}

The result (\ref{verlem1}) is an immediate consequence of Lemma
3.4.2 in \cite{vw96}. To show (\ref{verlem2}), we first apply the
symmetrization arguments used in the proof of Lemma \ref{bmodu}.
For sufficiently small $\delta$, the left-hand side of
(\ref{ratecon5}) is bounded by
%
%e90 ###
%
\begin{equation}\label{bound2}
2E_{\XW}^{o} \Biggl\|\frac{1}{\sqrt{n}}\sum
_{i=1}^{n}W_{ni}Y_{ni} \Biggr\|
_{\mathcal{V}_{\delta}},
\end{equation}
where $W_{ni}$'s are the assumed bootstrap weights and
\[
Y_{ni}=\bigl\{(\delta_{X_{i}}-P_X)\bigl(v(\theta,\eta)-v(\theta,\eta
_{n})\bigr)\dvtx d(\eta
,\eta_{n})\leq\delta,
\|\theta-\theta_{0}\|\leq\delta\bigr\}.
\]
Next, the multiplier inequality (\ref{mulineq}) is employed to
further bound (\ref{bound2}). In view of~(\ref{mulineq}), we
need only to figure out the upper bound for
%
%e91 ###
%
\begin{equation}\label{uppbou1}
E_X^{o}\|Y_{n1}\|_{\mathcal{V}_{\delta}}
\end{equation}
and
%
%e92 ###
%
\begin{equation}\label{uppbou2}
\max_{n_{0}\leq i\leq
n}E^{o}_X \Biggl\|\frac{1}{\sqrt{i}}\sum_{j_{0}+1}^{i}Y_{nj} \Biggr\|
_{\mathcal{V}_{\delta}}
\end{equation}
for some $n_{0}\geq1$ given assumptions W3 and W4 on the
bootstrap weights. By a similar argument as in the proof of
Lemma \ref{bmodu}, we know
\[
E^{o}_X\|Y_{n1}\|_{\mathcal{V}_{\delta}}\leq2E^{o}_XV_{n}(X_{1}),
\]
where $V_{n}$ is the envelop function of the class
$\mathcal{V}_{\delta}$ defined in (\ref{addass2}). The assumption
(\ref{inter15}), together with the analysis of assumption SB1,
implies that
$\lim\sup_{n}E^{o}_X\|Y_{n1}\|_{\mathcal{V}_{\delta}}<\infty$.
Next, Lemma 3.4.2 in \cite{vw96} implies that
\[
E^{o}_X\|\mathbb{G}_{n}\|_{\mathcal{V}_{\delta}}\leq
K(\delta,\mathcal{V}_{\delta},L_{2}(P)) \biggl(1+\frac{K(\delta
,\mathcal
{V}_{\delta},L_{2}(P))}{\delta^{2}\sqrt{n}} \biggr).
\]
By the triangular inequality, we know that (\ref{uppbou2}) has the
same upper bound as
$E^{o}_X\|\mathbb{G}_{n}\|_{\mathcal{V}_{\delta}}$. This concludes
the proof of (\ref{verlem2}).
\end{appendix}

\section*{Acknowledgments}
The authors thank Professor Anirban DasGupta for continuous
encouragement and Professors Michael Kosorok and Jon Wellner for many
helpful discussions. The authors also thank the Co-editor Susan Murphy
and two referees for insightful comments which led to important
improvements over an earlier draft.

\printaddresses

\end{document}